\documentclass[twocolumn]{autart}%
\usepackage{color}
\usepackage{graphicx}
\usepackage{epsfig}
\usepackage{cite}
\usepackage{amsmath}
\usepackage{amssymb}
\usepackage{mathrsfs}
\usepackage{amsfonts}
\usepackage{dsfont}
\usepackage{graphicx}
\usepackage{mathrsfs}
\usepackage[round]{natbib}%
\providecommand{\U}[1]{\protect\rule{.1in}{.1in}}
\providecommand{\U}[1]{\protect\rule{.1in}{.1in}}

\newtheorem{innercustomthm}{Theorem}
\newenvironment{customthm}[1]
  {\renewcommand\theinnercustomthm{#1}\innercustomthm}
  {\endinnercustomthm}
\newtheorem{innercustomlem}{Lemma}
\newenvironment{customlem}[1]
  {\renewcommand\theinnercustomlem{#1}\innercustomlem}
  {\endinnercustomlem}
\newtheorem{innercustomdef}{Definition}
\newenvironment{customdef}[1]
  {\renewcommand\theinnercustomdef{#1}\innercustomdef}
  {\endinnercustomdef}
\newtheorem{innercustomass}{Assumption}
\newenvironment{customass}[1]
  {\renewcommand\theinnercustomass{#1}\innercustomass}
  {\endinnercustomass}
\newtheorem{innercustomcor}{Corollary}

\newtheorem{innercustomrem}{Remark}
\newenvironment{customrem}[1]
  {\renewcommand\theinnercustomrem{#1}\innercustomrem}
  {\endinnercustomcor}
\allowdisplaybreaks[0]
\begin{document}
\begin{frontmatter}
\title{A Decentralized Energy-Optimal Control Framework for
Connected Automated Vehicles at Signal-Free Intersections\thanksref{footnoteinfo}} 
\thanks[footnoteinfo]{This
research was supported by the US Department of Energy's (DOE) SMART Mobility
Initiative. The work of Cassandras and Zhang is supported in part by NSF under
grants CNS-1239021, ECCS-1509084, CNS-1645681, and IIP-1430145, by AFOSR under
grant FA9550-15-1-0471, by The MathWorks and by Bosch.}
\author[Paestum]{Andreas A. Malikopoulos}\ead{andreas@udel.edu},    
\author[Rome]{Christos G. Cassandras}\ead{cgc@bu.edu},               
\author[Rome]{Yue J. Zhang}\ead{joycez@bu.edu}  
\address[Paestum]{Department of Mechanical Engineering, University of Delaware, 126 Spencer Lab, 130 Academy Street, Newark, DE, 19716, USA}  
\address[Rome]{Division of Systems Engineering and Center for Information and Systems Engineering, Boston University, 15 Saint Mary's Street,
Brookline, MA, 02446, USA}             
\begin{abstract}                          
We address the problem of optimally controlling connected and automated vehicles
(CAVs) crossing an urban intersection without any explicit traffic signaling,
so as to minimize energy consumption subject to a throughput maximization
requirement. We show that the solution of the throughput maximization problem
depends only on the hard safety constraints imposed on CAVs and its
structure enables a decentralized optimal control problem
formulation for energy minimization. We present a complete analytical solution
of these decentralized problems and derive conditions under which feasible
solutions satisfying all safety constraints always exist. The effectiveness of
the proposed solution is illustrated through simulation which shows
substantial dual benefits of the proposed decentralized framework by allowing
CAVs to conserve momentum and fuel while also improving travel time.
\end{abstract}
\end{frontmatter}

\section{Introduction}

\label{sec:1} Next generation transportation networks are typical
cyber-physical systems where event-driven components monitor and control
physical entities online. We are currently witnessing an increasing
integration of energy, transportation, and cyber networks, which, coupled with
human interactions, is giving rise to a new level of complexity in the
transportation network and necessitates new control and optimization approaches.

The alarming state of current transportation systems is well documented. In
2014, congestion caused vehicles in urban areas to spend 6.9 billion
additional hours on the road at a cost of an extra 3.1 billion gallons of
fuel, resulting in a total cost estimated at \$160 billion; see
\citet{Schrank2015}. From a control and optimization standpoint, the challenge
is to develop mechanisms that expand capacity \emph{without} affecting the
existing road infrastructure, specifically by tighter spacing of vehicles in
roadways and better control at the weakest links of a transportation system:
the bottleneck points defined by intersections, merging roadways, and speed
reduction zones; see \citet{Malikopoulos2013}, \citet{Margiotta2011}. An
automated highway system (AHS) can alleviate congestion, reduce energy use and
emissions, and improve safety by significantly increasing traffic flow as a
result of closer packing of automatically controlled vehicles. Forming
\textquotedblleft platoons\textquotedblright\ of vehicles traveling at high
speed is a popular system-level approach to address traffic congestion that
gained momentum in the 1990s; see \citet{Shladover1991,Rajamani2000}. More
recently, a study in \citet{Ratti2016} indicated that transitioning from
intersections with traffic lights to autonomous ones has the potential of
doubling capacity and reducing delays.

Connected and automated vehicles (CAVs) provide the most intriguing
opportunity for enabling users to better monitor transportation network
conditions and to improve traffic flow. CAVs can be controlled at different
transportation segments, e.g., intersections, merging roadways, roundabouts,
speed reduction zones and can assist drivers in making better operating
decisions to improve safety and reduce pollution, energy consumption, and
travel delays. One of the very early efforts in this direction was proposed in
\citet{Athans1969} and \citet{Levine1966} where the merging problem was
formulated as a linear optimal regulator to control a single string of
vehicles. \citet{Varaiya1993} has also discussed extensively the key features
of an automated intelligent vehicle-highway system (IVHS) and proposed a
related control system architecture.

In this paper, we address the problem of optimally controlling CAVs crossing
an urban intersection without any explicit traffic signaling so as to minimize
energy consumption subject to a throughput maximization requirement and to
hard safety constraints. The implications of this approach are that vehicles
do not have to come to a full stop at the intersection, thereby conserving
momentum and fuel while also improving travel time. Moreover, by optimizing
each vehicle's acceleration/deceleration, we minimize transient engine
operation, thus we have additional benefits in fuel consumption. Several
research efforts have been reported in the literature proposing either
\emph{centralized} (if there is at least one task in the system that is
globally decided for all vehicles by a single central controller) or
\emph{decentralized} approaches for coordinating CAVs at intersections.
\citet{Dresner2004} proposed the use of a centralized reservation scheme to
control a single intersection of two roads with no turns allowed. Since then,
numerous centralized approaches have been reported in the literature, e.g.,
\citet{Dresner2008,DeLaFortelle2010, Huang2012}, to achieve safe and efficient
control of traffic through intersections. Some approaches have focused on
coordinating vehicles to improve the travel time, e.g.,
\citet{Zohdy2012,Yan2009,Zhu2015a}. Others have considered minimizing the
overlap in the position of vehicles inside the intersection rather than
arrival time; see \citet{Lee2013}. \citet{Kim2014} proposed an approach based
on model predictive control that allows each vehicle to optimize its
movement locally with respect to any objective of interest.
\citet{Miculescu2014} used queueing theory and modeled the problem as a
polling system that determines the sequence of times assigned to the vehicles
on each road.

In decentralized approaches, each vehicle determines its own control policy
based on the information received from other vehicles on the road or from a
coordinator. \citet{Alonso2011} proposed two conflict resolution schemes in
which an autonomous vehicle can make a decision about the appropriate order of
crossing the intersection to avoid collision with other manually driven
vehicles. \citet{Colombo2014} constructed the invariant set for the control
inputs that ensure lateral collision avoidance. A detailed discussion of
research efforts in this area can be found in \citet{Malikopoulos2016a}.

The first contribution of the paper is the formulation of an energy minimization optimal control
problem for CAVs where the time for each CAV to cross the intersection is
first determined as the solution of a throughput maximization problem. We show that
the solution structure of the latter problem enables a decentralized energy
minimization optimal control problem formulation whose terminal time depends
only on a \textquotedblleft neighboring\textquotedblright\ CAV set. An
analytical solution of each CAV's optimal control problem without considering
state and control constraints was presented in \citet{Rios-Torres2015},
\citet{Rios-Torres2}, \citet{Ntousakis:2016aa} for CAVs at highway on-ramps,
and in \citet{ZhangMalikopoulosCassandras2016} for two adjacent intersections.
Unlike all these prior formulations, we specify the explicit connection
between the energy minimization and throughput maximization problems, do not
impose constraints on the terminal CAV speeds, and present a complete
analytical solution that includes all state and control constraints. Ensuring
that a \emph{feasible} solution to each CAV decentralized optimal control
problem exists is nontrivial, as discussed in \citet{Zhang2016}. Thus, another
contribution is showing that this solution depends on the arrival time of a
CAV at a \textquotedblleft control zone\textquotedblright\ defined for the
intersection and on its initial speed and then providing a proof (not given in
\citet{Zhang2016}) of the existence of a nonempty feasibility region in the
space defined by this arrival time and initial speed.

The paper is organized as follows. In Section II, we introduce the modeling
framework, formulate the energy-minimization optimal control problem and
establish its connection to throughput maximization. In Section III, we
present the decentralized control framework, derive a closed-form analytical
solution for each decentralized problem, and show the existence of feasible
solutions ensuring that all safety constraints remain inactive. Finally, we
provide simulation results in Section IV illustrating the effectiveness of the
proposed solution in terms of significant reductions in both fuel consumption
and travel time. Concluding remarks are given in Section V.


\section{Problem Formulation}

\label{sec:2}

\begin{figure}[ptb]
\centering
\includegraphics[width=3.4 in]{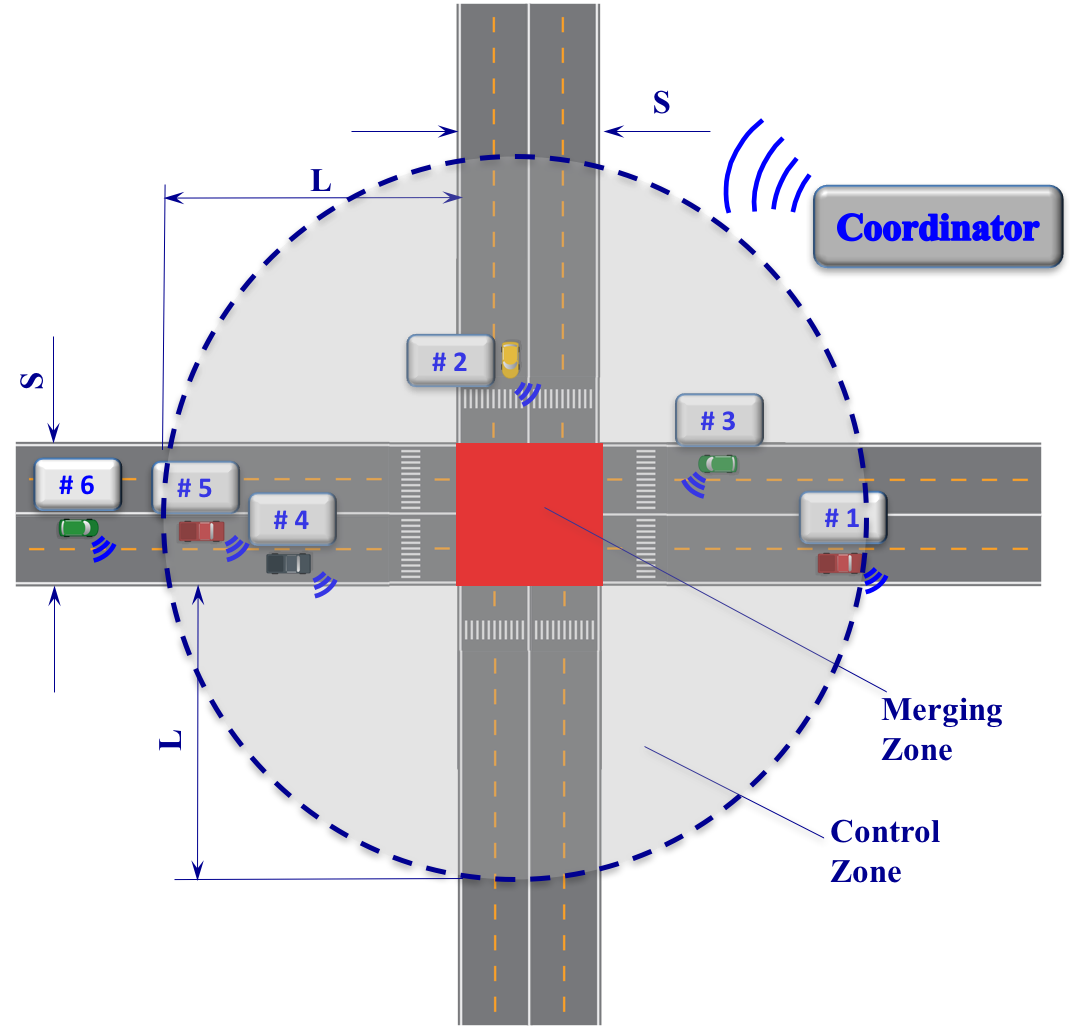} \caption{Intersection with
connected and automated vehicles.}%
\label{fig:1}%
\end{figure}

We consider an intersection (Fig. \ref{fig:1}) where the region at its center
is called \textit{Merging Zone} (MZ) and is the area of potential lateral
collision of vehicles. Although this is not restrictive, we consider the MZ to
be a square of side $S$. The intersection has a \textit{Control Zone} (CZ) and
a coordinator that can communicate with the vehicles traveling inside the CZ.
Note that the coordinator is not involved in any decision for any CAV and only
enables communication of appropriate information among CAVs. The distance from
the entry of the CZ to the entry of the MZ is $L$ and it is assumed to be the
same for all CZ entry points. The value of $L$ depends on the coordinator's
communication range capability with the CAVs, while $S$ is the physical length
of a typical intersection. In this paper, we limit ourselves to the case of no
lane changes and no turns allowed.

Let $N(t)\in\mathbb{N}$ be the number of CAVs inside the CZ at time
$t\in\mathbb{R}^{+}$ and $\mathcal{N}(t)=\{1,\ldots,N(t)\}$ be a queue which
designates the order in which these vehicles will be entering the MZ. Thus,
letting $t_{i}^{m}$ be the assigned time for vehicle $i$ to enter the MZ, we
require that
\begin{equation}
t_{i}^{m}\geq t_{i-1}^{m},~\forall i\in\mathcal{N}(t),~i>1. \label{eq:fifo}%
\end{equation}
There is a number of ways to satisfy \eqref{eq:fifo}. For example, we may
impose a strict first-in-first-out queueing structure, where each vehicle must
enter the MZ in the same order it entered the CZ. More generally, however,
$t_{i}^{m}$ may be determined for each vehicle $i$ at time $t_{i}^{0}$ when
the vehicle enters the CZ and $\mathcal{N}(t_{i}^{0})=\{1,\ldots,i-1\}$. If
$t_{i}^{m}>t_{i-1}^{m}$, then the order in the queue is preserved. If, on the
other hand, there exists some $j\in\mathcal{N}(t_{i}^{0})$, where $j<i-1$,
such that $t_{j}^{m}>t_{i}^{m}>t_{j-1}^{m}$, then the order is updated so that
CAV $i$ is placed in the $j$th queue position. The policy through which the
order (\textquotedblleft schedule\textquotedblright) is specified may be the
result of a higher level optimization problem as long as the condition
$t_{i}^{m}\geq t_{i-1}^{m}$ is preserved in between CAV arrival events at the
CZ. In what follows, we will adopt a specific scheme for determining
$t_{i}^{m}$ (upon arrival of CAV $i$) based on our problem formulation,
without affecting $t_{1}^{m},\ldots,t_{i-1}^{m}$, but we emphasize that our
analysis is not restricted by the policy designating the order of the vehicles
within the queue $\mathcal{N}(t)$.

\subsection{Vehicle Model, Constraints, and Assumptions}

For simplicity, we represent the dynamics of each CAV $i\in\mathcal{N}(t)$,
moving along a specified lane through second order dynamics%
\begin{equation}%
\begin{split}
\dot{p}_{i}  &  =v_{i}(t),\text{ \ \ }p_{i}(t_{i}^{0})=0\\
\dot{v}_{i}  &  =u_{i}(t),\text{ \ \ }v_{i}(t_{i}^{0})\text{ given}%
\end{split}
\label{eq:model2}%
\end{equation}
where $t_{i}^{0}$ is the time when CAV $i$ enters the CZ, and $p_{i}%
(t)\in\mathcal{P}_{i}$, $v_{i}(t)\in\mathcal{V}_{i}$, $u_{i}(t)\in
\mathcal{U}_{i}$ denote the position, speed and acceleration/deceleration
(control input) of each CAV $i$ inside the CZ. The sets $\mathcal{P}_{i}$,
$\mathcal{V}_{i}$ and $\mathcal{U}_{i}$, $i\in\mathcal{N}(t),$ are complete
and totally bounded subsets of $\mathbb{R}$. The state space $\mathcal{P}%
_{i}\times\mathcal{V}_{i}$ is closed with respect to the induced topology,
thus, it is compact.

We need to ensure that for any initial time and state $(t_{i}^{0},p_{i}%
^{0},v_{i}^{0})$ and every admissible control $u(t)$, the system
(\ref{eq:model2}) has a unique solution $(p_{i}(t),v_{i}(t))$ on some interval
$[t_{i}^{0},t_{i}^{m}]$, where $t_{i}^{m}$ is the time that vehicle
$i\in\mathcal{N}(t)$ enters the MZ. To ensure that the control input and
vehicle speed are within a given admissible range, the following constraints
are imposed:
\begin{equation}%
\begin{split}
u_{i,min}  &  \leqslant u_{i}(t)\leqslant u_{i,max},\quad\text{and}\\
0  &  \leqslant v_{min}\leqslant v_{i}(t)\leqslant v_{max},\quad\forall
t\in\lbrack t_{i}^{0},t_{i}^{m}],
\end{split}
\label{speed_accel constraints}%
\end{equation}
where $u_{i,min}$, $u_{i,max}$ are the minimum and maximum control inputs
(maximum deceleration/acceleration) for each vehicle $i\in\mathcal{N}(t)$, and
$v_{min}$, $v_{max}$ are the minimum and maximum speed limits respectively.
For simplicity, in the sequel we do not consider vehicle diversity and thus
set $u_{i,min}=u_{min}$ and $u_{i,max}=u_{max}$.

\begin{customdef}{1}
Depending on its physical location inside the CZ, CAV $i-1\in\mathcal{N}(t)$
belongs to only one of the following four subsets of $\mathcal{N}(t)$ with
respect to CAV $i$: 1) $\mathcal{R}_{i}(t)$ contains all CAVs traveling on the
same road as $i$ and towards the same direction but on different lanes (e.g.,
$\mathcal{R}_{6}(t)$ contains CAV 4 in Fig. \ref{fig:1}), 2) $\mathcal{L}%
_{i}(t)$ contains all CAVs traveling on the same road and lane as vehicle $i$
(e.g., $\mathcal{L}_{6}(t)$ contains CAV 5 in Fig. \ref{fig:1}), 3)
$\mathcal{C}_{i}(t)$ contains all CAVs traveling on different roads from $i$
and having destinations that can cause collision at the MZ, (e.g.,
$\mathcal{C}_{6}(t)$ contains CAV 2 in Fig. \ref{fig:1}), and 4)
$\mathcal{O}_{i}(t)$ contains all CAVs traveling on the same road as $i$ and
opposite destinations that cannot, however, cause collision at the MZ (e.g.,
$\mathcal{O}_{6}(t)$ contains CAV 3 in Fig. \ref{fig:1}). \label{def:2}
\end{customdef}

Based on this definition, it is clear that a rear-end collision can only arise
if CAV $k\in\mathcal{L}_{i}(t)$ is directly ahead of $i$. Thus, to ensure the
absence of any rear-end collision, we assume a predefined safe distance
$\delta<S$ and impose the rear-end safety constraint
\begin{equation}
s_{i}(t)=p_{k}(t)-p_{i}(t)\geqslant\delta,~\forall t\in\lbrack t_{i}^{0}%
,t_{i}^{f}]\text{, \ }k\in\mathcal{L}_{i}(t) \label{eq:rearend}%
\end{equation}
where $t_{i}^{f}$ is the time that CAV $i\in\mathcal{N}(t)$ exits the MZ. 
The rear-end safety constraint is usually expressed in terms of the allowable headway [\cite{Rajamani2012}], i.e., a time gap that is a function of speed. However, since we consider urban intersections, the average speed does not exhibit significant variations. Therefore, we can translate the allowable headway to a safe inter-vehicle distance. In
the rest of the paper, we reserve the symbol $k$ to denote the CAV which is
physically immediately ahead of $i$ in the same lane.

A lateral collision involving CAV $i$ may occur only if some CAV $j\neq i$
belongs to $\mathcal{C}_{i}(t)$. This leads to the following definition:


\begin{customdef}{2}
For each CAV $i\in\mathcal{N}(t)$, we define the set $\Gamma_{i}$ that
includes all time instants when a lateral collision involving CAV $i$ is
possible:
\begin{equation}
\Gamma_{i}\triangleq\Big\{t~|~t\in\lbrack t_{i}^{m},t_{i}^{f}]\Big\}.
\end{equation}

\end{customdef}

Consequently, to avoid a lateral collision for any two vehicles $i,j\in
\mathcal{N}(t)$ on different roads, the following constraint should hold
\begin{equation}
\Gamma_{i}\cap\Gamma_{j}=\varnothing,\text{ \ \ \ }\forall t\in\lbrack
t_{i}^{m},t_{i}^{f}]\text{, \ }j\in\mathcal{C}_{i}(t). \label{eq:lateral}%
\end{equation}
This constraint implies that no two CAVs from different roads which may lead
to a lateral collision are allowed to be in the MZ at the same time. If the
length of the MZ is large, then this constraint might not be realistic, but it
can be modified appropriately as described in Remark \ref{remark:1}.

In the modeling framework described above, we impose the following assumptions:


\begin{customass}{1}
For CAV $i$, none of the constraints (\ref{speed_accel constraints}%
)-(\ref{eq:rearend}) is active at $t_{i}^{0}$. \label{ass:feas}
\end{customass}


\begin{customass}{2}
The speed of the CAVs inside the MZ is constant, i.e., $v_{i}(t)=v_{i}%
(t_{i}^{m})=v_{i}(t_{i}^{f})$, $\ \forall t\in\lbrack t_{i}^{m},t_{i}^{f}]$
\label{ass:4} This implies that
\begin{equation}
t_{i}^{f}=t_{i}^{m}+\frac{S}{v_{i}(t_{i}^{m})}. \label{eq:time}%
\end{equation}

\end{customass}


\begin{customass}{3}
Each CAV $i$ has proximity sensors and can measure local information without
errors or delays. \label{ass:sensor}
\end{customass}

Assumption \ref{ass:feas} ensures that the initial state and control input are
feasible. Enforcing this is nontrivial and we address the issue in Section
\ref{sec:4c}. The second assumption is intended to enhance safety awareness,
but it could be modified appropriately, if necessary, as discussed in Section
\ref{sec:2b}. The third assumption may be strong, but it is relatively
straightforward to relax as long as the noise in the measurements and/or
delays is bounded. For example, we can determine upper bounds on the state
uncertainties as a result of sensing or communication errors and delays, and
incorporate these into more conservative safety constraints.

For simplicity of notation in the remainder of the paper, we will write
$v_{i}(t_{i}^{0})\equiv v_{i}^{0}$, $v_{i}(t_{i}^{m})\equiv v_{i}^{m}$ and
$v_{i}(t_{i}^{f})\equiv v_{i}^{f}$.


\subsection{Energy Minimization Problem Formulation}

\label{sec:2b} We begin by considering the controllable
acceleration/deceleration $u_{i}(t)$ of each CAV $i$ which minimizes the
following cost functional:
\begin{gather}
J_{i}(u_{i}(t),t_{i}^{m})=\int_{t_{i}^{0}}^{t_{i}^{m}}C_{i}%
(u_{i}(t))dt,\label{eq:functional}\\
\text{subject to}%
:\eqref{eq:model2},\eqref{speed_accel constraints},\eqref{eq:rearend},\eqref{eq:lateral},\text{
}p_{i}(t_{i}^{0})=0\text{, }p_{i}(t_{i}^{m})=L,\nonumber\\
\text{and given }t_{i}^{0}\text{, }v_{i}^{0}\text{, }t_{i}^{m}.\nonumber
\end{gather}
We view $C_{i}(u_{i}(t))$ as a measure of the energy, which is a function of
the control input (acceleration/deceleration) consumed by CAV $i$ in traveling
between $p_{i}(t_{i}^{0})=0$ and $p_{i}(t_{i}^{m})=L$; see
\citet{Malikopoulos2011}. A special case arises when the cost function is the
$L^{2}$-norm of the control input in $[t_{i}^{0},t_{i}^{m}]$ and $C_{i}%
(u_{i}(t))=\frac{1}{2}u_{i}^{2}(t)$. In this case, we minimize transient
engine operation, thus we can have direct benefits in fuel consumption and
emissions since internal combustion engines are optimized over steady state
operating points (constant torque and speed); see \citet{Rios-Torres2} and
\citet{malikopoulos2013stochastic}. In this problem, $t_{i}^{0}$, $v_{i}^{0}$
are known upon arrival of CAV $i$ at the CZ and $t_{i}^{m}$ is also specified.
Clearly, not all $t_{i}^{m}$ can satisfy the safety constraints
\eqref{eq:rearend}~and~\eqref{eq:lateral}. Moreover, in general, a value of
$t_{i}^{m}$ that satisfies \eqref{eq:rearend}~and~\eqref{eq:lateral} may
depend on other CAVs $j\neq i$; therefore, it may not be possible for CAV $i$
to solve (\ref{eq:functional}) in a \emph{decentralized} manner, i.e., based
only on local information. We address the question of specifying appropriate
$t_{i}^{m}$ for each instance of (\ref{eq:functional}) in what follows.

Before proceeding, we note that the obvious unconstrained solution to
(\ref{eq:functional}) is $u_{i}^{\ast}(t)=0$ for all $t\in\lbrack t_{i}%
^{0},t_{i}^{m}]$. This applies to $i=1$ since, in this case,
\eqref{eq:rearend}~and~\eqref{eq:lateral} are inactive, since it is not
constrained by any prior CAV in the queue, and $t_{1}^{m}$ variable. This
also implies that $v_{1}^{\ast}(t)=v_{i}^{0}$ for all $t\in\lbrack t_{i}%
^{0},t_{i}^{m}]$ and $t_{1}^{m}=L/v_{i}^{0}$.

We now turn our attention to the problem of maximizing the traffic throughput
at the intersection, in terms of minimizing the gaps between the vehicles in a
given queue $\mathcal{N}(t)$ (see Fig. \ref{fig:1}), under the hard safety
constraints (\ref{eq:rearend}) and (\ref{eq:lateral}). Thus, setting
$\mathbf{t}_{(2:N(t))}=[t_{2}^{m}\ldots t_{N(t)}^{m}]$, we define the
following optimization problem:
\begin{gather}
\min_{\mathbf{t}_{(2:N(t))}}\sum_{i=2}^{N(t)}\Big(t_{i}^{m}-t_{i-1}%
^{m}\Big)=\min_{\mathbf{t}_{N(t)}}\Big(t_{N(t)}^{m}-t_{1}^{m}%
\Big),\label{eq:central}\\
\text{subject to}%
:\eqref{eq:fifo},\eqref{speed_accel constraints},\eqref{eq:rearend},\eqref{eq:lateral}.\nonumber
\end{gather}
where $t_{1}^{m}$ is not included since it is obtained from the solution of
(\ref{eq:functional}) when $i=1$, i.e., $t_{1}^{m}=L/v_{i}^{0}$. The
equivalence between the two expressions in \eqref{eq:central} (due to the
cancellation of all terms in the sum except the first and last) reflects the
equivalence between minimizing the total time to process all CAVs in the queue
and the average interarrival time of CAVs at the MZ.

As stated in (\ref{eq:central}), the problem does not incorporate constraints
on $t_{i}^{m}$, $i=2,\ldots,N(t)$, that are imposed by the CAV dynamics. In
other words, we should write $t_{i}^{m}=t_{i}^{m}(\mathbf{u}_{(1:i)}(t))$
where $\mathbf{u}_{(1:i)}(t)=[u_{1}(t;t_{1}^{m})\ldots u_{i}(t;t_{i}^{m})]$
denotes the controls applied to all CAVs $i=1,\ldots,N(t)$ over $[t_{i}%
^{0},t_{i}^{m}]$ for any given $t_{i}^{0},t_{i}^{m}$. Let $\mathscr{A}_{i}$
denote a set of feasible controls:%
\begin{gather}
\mathscr{A}_{i}\triangleq\Big\{u_{i}(t;t_{i}^{m})\in\mathcal{U}_{i}%
~\text{subject to:}\label{FeasibleControlSet}\\
\eqref{eq:fifo},\eqref{eq:model2},\eqref{speed_accel constraints},\eqref{eq:rearend},\eqref{eq:lateral},~p_{i}%
(t_{i}^{0})=0\text{, }p_{i}(t_{i}^{m})=L,\nonumber\\
\text{and given }t_{i}^{0}\text{, }v_{i}^{0}\text{, }t_{i}^{m}\Big\}.\nonumber
\end{gather}
Then, we rewrite (\ref{eq:central}) as
\begin{gather}
\min_{\mathbf{t}_{(2:N(t))}}\sum_{i=2}^{N(t)}\Big(t_{i}^{m}(\mathbf{u}%
_{(1:i)}(t))-t_{i-1}^{m}(\mathbf{u}_{(1:i-1)}(t))\Big)\label{eq:central1}\\
=\min_{\mathbf{t}_{N(t)}}\Big(t_{N(t)}^{m}(\mathbf{u}_{(1:N(t))}(t))-t_{1}%
^{m}(\mathbf{u}_{(1)}(t))\Big),\nonumber
\end{gather}%
\[
\text{subject to}:u_{i}(t;t_{i}^{m})\in\mathscr{A}_{i},~\forall i\in
\mathcal{N}(t),~\eqref{eq:fifo},\text{ }\eqref{speed_accel constraints},\text{
}\eqref{eq:rearend},~\eqref{eq:lateral}.
\]

\begin{customrem}{1}
As pointed out earlier, the solution of (\ref{eq:functional}) for $i=1$ is
$u_{1}^{\ast}(t)=0$ resulting in $v_{1}^{\ast}(t)=v_{i}^{0}$ and $t_{1}%
^{m\ast}=L/v_{i}^{0}$. On the other hand, if we were to solve
(\ref{eq:central1}) for $i=1$ setting $t_{0}^{m}=0$, the solution would be
$t_{1}^{m\ast}=t_{1}^{m}=L/v_{max}$. This indicates a degree of freedom in the
selection of $t_{1}^{m}$ which can be used to trade off the energy
minimization and throughput maximization (congestion reduction) objectives.
Thus, $t_{1}^{m}$ may be viewed as a parameter one can adjust to solve the
subsequent CAV problems placing a desired amount of emphasis on throughput
relative to energy consumption. \label{Remark_t1choice}
\end{customrem}

The solution of (\ref{eq:central1}) provides a sequence $\{t_{2}^{m\ast
},\ldots,t_{N(t)}^{m\ast}\}$ which designates the MZ arrival times of all CAVs
in the current queue so as to minimize the total time needed for them to clear
the intersection (recalling Assumption 2, the time through the MZ is fixed),
hence maximizing the throughput over the current $N(t)$ CAVs. This solution
may then be used in (\ref{eq:functional}) to specify the terminal time of each
energy minimization problem. In what follows, we show that this solution has a
simple iterative structure and depends only on the hard safety constraints
\eqref{eq:rearend}~and~\eqref{eq:lateral}, as well as the state and control
constraints \eqref{speed_accel constraints}. We begin by ignoring the latter
to obtain the following result.

\begin{customlem}{1}
Suppose that the constraints \eqref{speed_accel constraints} are inactive in
(\ref{eq:central1}). Then, the solution $\mathbf{t}^{\ast}=[t_{2}^{m^{\ast}%
},\ldots,t_{N}^{m^{\ast}}]$ of problem (\ref{eq:central1}) is determined
through the following recursive structure over $i=2,\ldots,N$: \label{lemma 2}%
\begin{equation}
t_{i}^{m^{\ast}}=\left\{
\begin{array}
[c]{ll}%
\max\{t_{i-1}^{m^{\ast}},t_{k}^{m^{\ast}}+\frac{\delta}{v_{k}^{m}}\} &
\text{if }i-1\in\mathcal{R}_{i}(t)\cup\mathcal{O}_{i}(t)\\
t_{i-1}^{m^{\ast}}+\frac{\delta}{v_{i-1}^{m}} &
\mbox{if $i-1\in\mathcal{L}_{i}$}\\
t_{i-1}^{m^{\ast}}+\frac{S}{v_{i-1}^{m}} & \mbox{if $i-1\in\mathcal{C}_{i}$}
\end{array}
\right.  \label{Lemma}%
\end{equation}
where $k=\max\{j:j\in\mathcal{L}_{i}(t),$ $j=1,\ldots,i-2\}<i$ is the CAV
which is physically immediately ahead of $i$ in the same lane.
\end{customlem}

\textbf{Proof:} See Appendix.

\begin{customrem}{2}
The lateral collision constraint \eqref{eq:lateral} allows only one CAV at a
time to be inside the MZ. If the length of the MZ is large, however, then this
constraint may become overly conservative, since it results in dissipating
space and road capacity. The constraint can be modified appropriately and
(\ref{lemmacase3}) in \emph{Case 3} above can be rewritten as
\begin{equation}
t_{i}^{m^{\ast}}=t_{i-1}^{m^{\ast}}+\frac{r}{v_{i-1}^{m}} \label{eq:rem2}%
\end{equation}
with any desired distance $r<S$ between CAVs inside the MZ. \label{remark:1}
\end{customrem}

Next, we relax the assumption made in Lemma \ref{lemma 2} that constraints
\eqref{speed_accel constraints} are inactive in (\ref{eq:central1}) and derive
a recursive equation for the determination of $\mathbf{t}^{\ast}%
=[t_{2}^{m^{\ast}},\ldots,t_{N}^{m^{\ast}}]$.

\begin{customthm}{1}
The solution $\mathbf{t}^{\ast}=[t_{1}^{m^{\ast}},\ldots,t_{N}^{m^{\ast}}]$ of
problem (\ref{eq:central1}) is recursively determined through%
\begin{equation}
t_{i}^{m^{\ast}}=\left\{
\begin{array}
[c]{ll}%
t_{1}^{m^{\ast}} & \mbox{if $i=1$}\\
\text{max }\{t_{i-1}^{m^{\ast}},t_{k}^{m^{\ast}}+\frac{\delta}{v_{k}^{m}%
},t_{i}^{c}\} & \text{if }i-1\in\mathcal{R}_{i}(t)\cup\mathcal{O}_{i}(t)\\
\text{max }\{t_{i-1}^{m^{\ast}}+\frac{\delta}{v_{i-1}^{m}},t_{i}^{c}\} &
\mbox{if $i-1\in\mathcal{L}_{i}$}\\
\text{max }\{t_{i-1}^{m^{\ast}}+\frac{S}{v_{i-1}^{m}},t_{i}^{c}\} &
\mbox{if $i-1\in\mathcal{C}_{i}$}
\end{array}
\right.  \label{def:tf}%
\end{equation}
where $t_{i}^{c}=t_{i}^{1}\mathds{1}_{v_{i}^{m}=v_{max}}+t_{i}^{2}%
(1-\mathds{1}_{v_{i}^{m}=v_{max}})$ and
\begin{align}
t_{i}^{1}  &  =t_{i}^{0}+\frac{L}{v_{max}}+\frac{(v_{max}-v_{i}^{0})^{2}%
}{2u_{i,max}v_{max}}\label{ti1def}\\
t_{i}^{2}  &  =t_{i}^{0}+\frac{[2Lu_{i,max}+(v_{i}^{0})^{2}]^{1/2}-v_{i}^{0}%
}{u_{i,max}} \label{ti2def}%
\end{align}
\label{theo:1}
\end{customthm}

\textbf{Proof:} See Appendix.

It follows from Theorem \ref{theo:1} that $t_{i}^{m^{\ast}}$ is always
recursively determined from $t_{i-1}^{m^{\ast}}$ and $v_{i-1}^{m}$ and
possibly $t_{k}^{m^{\ast}},$ $v_{k}^{m}$ where $v_{i-1}^{m}$ and $v_{k}^{m}$
depend on the specific controls used when solving problem (\ref{eq:central1}).
However, note that there is no guarantee that there exist feasible controls
satisfying all constraints in (\ref{FeasibleControlSet}) over all $t\in\lbrack
t_{i}^{0},t_{i}^{m}]$. In fact, as we will discuss in Section \ref{sec:4c}, it
is easy to see that the safety constraint \eqref{eq:rearend} may not hold
depending on the initial conditions $(t_{i}^{0},v_{i}^{0})$ for CAV $i$. We
will show, however, in Theorem \ref{theo:3} that there exists a nonempty
feasible region $\mathcal{F}_{i}\subset R^{2}$ of initial conditions
$(t_{i}^{0},v_{i}^{0})$ such that $s_{i}(t)\geqslant\delta$ for all
$t\in(t_{i}^{0},t_{i}^{m})$ so that all safety constraints are guaranteed to
hold throughout $[t_{i}^{0},t_{i}^{m}]$.

We are now in a position to return to the energy minimization problem
(\ref{eq:functional}) with the value of $t_{i}^{m}$ for any $i=1,\ldots,N(t)$
specified through (\ref{def:tf}) in a recursive manner. This allows us to
solve these problems in a decentralized manner as detailed in the next section.

\section{Decentralized Framework}

\label{sec:4}

The results in the previous section allow us to address the optimal control
problem \eqref{eq:functional} within a decentralized framework. However, to
establish this framework, we need a communication structure
between CAVs with a \textquotedblleft coordinator\textquotedblright\ whose
task is to handle the information between them. In particular, when a CAV $i$
reaches the CZ of the intersection at some instant $t$, the coordinator
assigns to it a unique identity as follows.

Let $M(t)\in\mathbb{N}$ be the cumulative number of CAVs that have entered the
CZ by time $t$. Note that $M(t)$ is increasing in $t$ and can be reset to
$M(t)=0$ only if no CAVs are inside the CZ. The \textit{unique identity} that
the coordinator assigns to each CAV is a triplet $(w,i,j)$ where $w=M(t)+1$ is
a unique index, $i$ is the position of the vehicle in the current queue
$\mathcal{N}(t)$, and $j\in\{1,\ldots,4\}$ is an integer based on a one-to-one
mapping from $\{\mathcal{R}_{i}(t),$ $\mathcal{L}_{i}(t),$ $\mathcal{C}%
_{i}(t),$ $\mathcal{O}_{i}(t)\}$ onto $\{1,\ldots,4\}$ that indicates the
positional relationship between CAVs $i-1$ and $i$. If two or more CAVs enter
the CZ at the same time, then the coordinator assigns randomly the index $w$.


\begin{customdef}{3}
For each vehicle $i$ entering a CZ, we define the \textit{information set}
$Y_{i}(t)$ as
\begin{equation}
Y_{i}(t)\triangleq\Big\{p_{i}(t),v_{i}(t),w,\mathcal{Q}_{i},s_{i}%
(t),t_{i}^{m^{\ast}}\Big\},\forall t\in\lbrack t_{i}^{0},t_{i}^{m^{\ast}}],
\end{equation}
\label{infoset}where $p_{i}(t),v_{i}(t)$ are the position and speed of CAV $i$
inside the CZ; $w$ and $\mathcal{Q}_{i}\in\{1,\ldots,4\}$ are the unique index
and the queue subset (Definition \ref{def:2}) assigned to CAV $i$ by the
coordinator; and $s_{i}(t)=p_{k}(t)-p_{i}(t)$ is the distance between CAV $i$
and some CAV $k$ which is immediately ahead of $i$ in the same lane (recall
that we reserve the symbol $k$ to denote such a CAV relative to $i$). The last
element above, $t_{i}^{m^{\ast}}$, is the time targeted for CAV $i$ to enter
the MZ and is given in (\ref{def:tf}) depending on the value of $\mathcal{Q}%
_{i}$.
\end{customdef}

Note that once CAV $i$ enters the CZ, then immediately all information in
$Y_{i}(t)$ becomes available to $i$: $p_{i}(t),v_{i}(t)$ are read from its
sensors; $\mathcal{Q}_{i}$ is assigned by the coordinator, as is the value of
$k$ based on which $s_{i}(t)$ is also evaluated; $t_{i}^{m^{\ast}}$ can also
be computed at that time based on the information the vehicle $i$ receives
from $i-1$. The recursion on $t_{i}^{m^{\ast}}$ is initialized whenever a
vehicle enters the CZ. In this case, $t_{1}^{m}$ can be externally assigned as
the desired exit time of this vehicle whose behavior is unconstrained (as
discussed in the previous section). Thus, the time $t_{1}^{m}$ is available
through $Y_{1}(t)$.

Since the coordinator is not involved in any control decision, from Theorem
\ref{theo:1} we can formulate $N(t)$ sequential decentralized tractable
problems of the form \eqref{eq:functional} that may be solved online. As
already discussed, a special case of \eqref{eq:functional} arises when the
cost function is the $L^{2}$-norm of the control input in $[t_{i}^{0}%
,t_{i}^{m^{\ast}}]$ which we shall henceforth consider. Thus, the
decentralized problem for each CAV $i$ is formulated as follows:
\begin{gather}
\min_{u_{i}(t)}\frac{1}{2}\int_{t_{i}^{0}}^{t_{i}^{m^{\ast}}}u_{i}%
^{2}(t)~dt,\label{eq:decentral}\\
\text{subject to}%
:\eqref{eq:model2},\eqref{speed_accel constraints},\eqref{def:tf},~p_{i}%
(t_{i}^{0})=0\text{, }p_{i}(t_{i}^{m^{\ast}})=L,\text{and}\nonumber\\
\text{given }t_{i}^{0}\text{, }v_{i}^{0}.\nonumber
\end{gather}
Observe that we have omitted the rear end safety constraint
\eqref{eq:rearend}~and~the lateral collision constraint \eqref{eq:lateral}.
The latter applies to the MZ and affects (\ref{eq:decentral}) only at
$t=t_{i}^{m^{\ast}}$ which is implicitly handled by the selection of
$t_{i}^{m^{\ast}}$ in (\ref{def:tf}). The former, on the other hand, must be
satisfied for all $t\in\lbrack t_{i}^{0},t_{i}^{m^{\ast}}]$, whereas
(\ref{def:tf}) only guarantees that it is satisfied at $t=t_{i}^{m^{\ast}}$.
It is omitted here because we will show that the solution of
(\ref{eq:decentral}) guarantees that this constraint indeed holds throughout
$[t_{i}^{0},t_{i}^{m^{\ast}}]$ under proper initial conditions $(t_{i}%
^{0},v_{i}^{0})$; note that the constraint also holds in $[t_{i}^{m^{\ast}%
},t_{i}^{f}]$ under Assumption \ref{ass:4}.


\subsection{Analytical solution of the decentralized optimal control problem}

\label{sec:4b}

For the analytical solution of \eqref{eq:decentral} and its online
implementation, we apply Hamiltonian analysis under Assumption \ref{ass:feas},
i.e., when the CAVs enter the CZ none of the constraints is active. We stress
that this is not in general true. For example, a CAV may enter the CZ with
speed higher than the speed limit. In this case, a solution of the optimal
control problem is infeasible. A feasibility analysis for CAVs to satisfy such
initial conditions is discussed in Section \ref{sec:4c} where we show that a
feasible region $\mathcal{F}_{i}\subset R^{2}$ of initial conditions
$(t_{i}^{0},v_{i}^{0})$ for CAV $i$ exists such that $s_{i}(t)\geqslant\delta$
for all $t\in\lbrack t_{i}^{0},t_{i}^{m}]$; a feasibility enforcement analysis
to ensure the existence of feasible and optimal solutions is given in \citet{Zhang2016}.

From \eqref{eq:decentral}, the state equations \eqref{eq:model2}, and the
control/state constraints \eqref{speed_accel constraints}, for each vehicle
$i\in\mathcal{N}(t)$ the Hamiltonian function with the state and control
constraints adjoined is
\begin{gather}
H_{i}\big(t,p(t),v(t),u(t)\big)=\frac{1}{2}u_{i}^{2}+\lambda_{i}^{p}\cdot
v_{i}+\lambda_{i}^{v}\cdot u_{i}\nonumber\\
+\mu_{i}^{a}\cdot(u_{i}-u_{max})+\mu_{i}^{b}\cdot(u_{min}-u_{i})+\mu_{i}%
^{c}\cdot(v_{i}-v_{max})\nonumber\\
+\mu_{i}^{d}\cdot(v_{min}-v_{i}), \label{eq:16b}%
\end{gather}
where $\lambda_{i}^{p}$ and $\lambda_{i}^{v}$ are the costates, and $\mu^{T}$
is a vector of Lagrange multipliers with
\begin{equation}
\mu_{i}^{a}=\left\{
\begin{array}
[c]{ll}%
>0, & \mbox{$u_{i}(t) - u_{max} =0$},\\
=0, & \mbox{$u_{i}(t) - u_{max} <0$},
\end{array}
\right.  \label{eq:17a}%
\end{equation}%
\begin{equation}
\mu_{i}^{b}=\left\{
\begin{array}
[c]{ll}%
>0, & \mbox{$u_{min} - u_{i}(t) =0$},\\
=0, & \mbox{$u_{min} - u_{i}(t)<0$},
\end{array}
\right.  \label{eq:17b}%
\end{equation}%
\begin{equation}
\mu_{i}^{c}=\left\{
\begin{array}
[c]{ll}%
>0, & \mbox{$v_{i}(t) - v_{max} =0$},\\
=0, & \mbox{$v_{i}(t) - v_{max}<0$},
\end{array}
\right.  \label{eq:17c}%
\end{equation}%
\begin{equation}
\mu_{i}^{d}=\left\{
\begin{array}
[c]{ll}%
>0, & \mbox{$v_{min} - v_{i}(t)=0$},\\
=0, & \mbox{$v_{min} - v_{i}(t)<0$}.
\end{array}
\right.  \label{eq:17d}%
\end{equation}

The Euler-Lagrange equations become
\begin{equation}
\dot{\lambda}_{i}^{p}=-\frac{\partial H_{i}}{\partial p_{i}}=0, \label{eq:EL1}%
\end{equation}
and
\begin{equation}
\dot{\lambda}_{i}^{v}=-\frac{\partial H_{i}}{\partial v_{i}}=\left\{
\begin{array}
[c]{ll}%
-\lambda_{i}^{p}, & \mbox{$v_{i}(t) - v_{max} <0$}~\text{and}\\
& \mbox{$v_{min} - v_{i}(t)<0$},\\
-\lambda_{i}^{p}-\mu_{i}^{c}, & \mbox{$v_{i}(t) - v_{max} =0$},\\
-\lambda_{i}^{p}+\mu_{i}^{d}, & \mbox{$v_{min} - v_{i}(t)=0$},
\end{array}
\right.  \label{eq:EL2}%
\end{equation}
with boundary conditions $p_{i}(t_{i}^{0})=0$, $p_{i}(t_{i}^{m})=L$,
$\lambda_{i}^{v}(t_{i}^{m})= 0$, given initial conditions $t_{i}^{0}$,
$v_{i}(t_{i}^{0})$, and $t_{i}^{m}$ specified by (\ref{def:tf}). The necessary
condition for optimality is
\begin{equation}
\frac{\partial H_{i}}{\partial u_{i}}=u_{i}+\lambda_{i}^{v}+\mu_{i}^{a}%
-\mu_{i}^{b}=0, \label{eq:KKT1}%
\end{equation}
To address this problem, constrained and unconstrained arcs need to be pieced
together to satisfy (\ref{eq:EL1}) through (\ref{eq:KKT1}). 
Based on our state and control constraints \eqref{speed_accel constraints} and boundary conditions, the optimal solution is the result of different combinations of the following possible arcs.


\textbf{1. Control and State Constraints not Active. }In this case, we have
$\mu_{i}^{a}=\mu_{i}^{b}=\mu_{i}^{c}=\mu_{i}^{d}=0.$ Applying \eqref{eq:KKT1},
the optimal control is given by
\begin{equation}
u_{i}+\lambda_{i}^{v}=0,\quad i\in\mathcal{N}(t).\label{eq:17}\\
\end{equation}
and the Euler-Lagrange equations yield (\ref{eq:EL1}) and
\begin{equation}
\dot{\lambda}_{i}^{v}=-\frac{\partial H_{i}}{\partial v_{i}}=-\lambda_{i}^{p}.
\label{eq:19}%
\end{equation}
From (\ref{eq:EL1}) we have $\lambda_{i}^{p}=a_{i}$ and \eqref{eq:19} implies
$\lambda_{i}^{v}=-(a_{i}t+b_{i})$, where $a_{i}$ and $b_{i}$ are integration
constants. Consequently, the optimal control input (acceleration/deceleration)
as a function of time is given by
\begin{equation}
u_{i}^{\ast}(t)=a_{i}t+b_{i}. \label{eq:20}%
\end{equation}
Substituting this equation into the vehicle dynamics \eqref{eq:model2} we can
find the optimal speed and position for each vehicle, namely
\begin{equation}
v_{i}^{\ast}(t)=\frac{1}{2}a_{i}t^{2}+b_{i}t+c_{i} \label{eq:21}%
\end{equation}%
\begin{equation}
p_{i}^{\ast}(t)=\frac{1}{6}a_{i}t^{3}+\frac{1}{2}b_{i}t^{2}+c_{i}t+d_{i},
\label{eq:22}%
\end{equation}
where $c_{i}$ and $d_{i}$ are integration constants. These fours constants
above can be computed by using the initial and final conditions in
(\ref{eq:decentral}). In particular, using \eqref{eq:21} with the initial
condition $v_{i}(t_{i}^{0})=v_{i}^{0}$, \eqref{eq:22} with the initial and
terminal conditions $p_{i}(t_{i}^{0})=0,$ $p_{i}(t_{i}^{m})=L$, and the
boundary condition of the costate $\lambda_{i}^{v}(t_{i}^{m})=-u_{i}(t_{i}%
^{m})=0$, we can form the system of four equations of the form $\mathbf{T}%
_{i}\mathbf{b}_{i}=\mathbf{q}_{i}$:%

\begin{equation}
\left[
\begin{array}
[c]{cccc}%
\frac{1}{6}(t_{i}^{0})^{3} & \frac{1}{2}(t_{i}^{0})^{2} & t_{i}^{0} & 1\\
\frac{1}{2}(t_{i}^{0})^{2} & t_{i}^{0} & 1 & 0\\
\frac{1}{6}(t_{i}^{m})^{3} & \frac{1}{2}(t_{i}^{m})^{2} & t_{i}^{m} & 1\\
-t_{i}^{m} & -1 & 0 & 0
\end{array}
\right]  .\left[
\begin{array}
[c]{c}%
a_{i}\\
b_{i}\\
c_{i}\\
d_{i}%
\end{array}
\right]  =\left[
\begin{array}
[c]{c}%
p_{i}(t_{i}^{0})\\
v_{i}(t_{i}^{0})\\
p_{i}(t_{i}^{m})\\
\lambda_{i}^{v}(t_{i}^{m})
\end{array}
\right]  \label{eq:23}%
\end{equation}
where $t_{i}^{m}$ is specified by \eqref{def:tf}. Note that since
\eqref{eq:23} can be computed online, the controller may re-evaluate the four
constants in the form $a_{i}(t,p_{i},v_{i}),b_{i}(t,p_{i},v_{i}),c_{i}%
(t,p_{i},v_{i})$, and $d_{i}(t,p_{i},v_{i})$ at any time $t>t_{i}^{0}$ to get
\begin{equation}
\mathbf{b}_{i}(t,p_{i}(t),v_{i}(t))=(\mathbf{T}_{i})^{-1}.\mathbf{q}%
_{i}(t,p_{i}(t),v_{i}(t)) \label{eq:24}%
\end{equation}
and update \eqref{eq:20} as follows
\begin{equation}
u_{i}^{\ast}(t,p_{i}(t),v_{i}(t))=a_{i}(t,p_{i}(t),v_{i}(t))t+b_{i}%
(t,p_{i}(t),v_{i}(t)). \label{eq:25}%
\end{equation}
Thus, feedback can be indirectly provided through the re-calculation of the
vector $\mathbf{b}_{i}(t,p_{i}(t),v_{i}(t))$ in \eqref{eq:24}.


\textbf{2. Control Constraint Active,} $u_{i}^{\ast}(t)=u_{max}$. Suppose that
at time $t=t_{1}$, (\ref{eq:20}) becomes
\begin{equation}
u_{i}^{\ast}(t)=u_{max},~\forall t\geq t_{1}. \label{eq:cas1d}%
\end{equation}
while $v_{min}<v_{i}(t)<v_{max}$. In this case, the Hamiltonian is continuous
at $t=t_{1}$ (entry point of the control constrained arc). Substituting the
last equation into the vehicle dynamics \eqref{eq:model2}, we can find the
optimal speed and position of each vehicle, namely
\begin{gather}
v_{i}^{\ast}(t)=u_{max}~t+f_{i},\label{eq:cas1e}\\
p_{i}^{\ast}(t)=\frac{1}{2}u_{max}~t^{2}+f_{i}~t+e_{i},~\forall t\geq t_{1}
\label{eq:cas1f}%
\end{gather}
where $f_{i}$ and $e_{i}$ are constants of integration that can be computed
easily since we know the speed and position of the vehicle at time $t=t_{1}$.


\textbf{3. Control and State Constraints Active, }$u_{i}(t)=u_{max}$ and
$v_{i}(t)=v_{max}$. Suppose that at time $t=t_{2}>t_{1}$ (exit point of the
control constrained arc and entry point of the state variable constrained
arc), \eqref{eq:cas1e} becomes $v_{i}^{\ast}(t)=v_{max}$. Then from
\eqref{eq:model2} we have $\dot{v}_{i}^{\ast}=u_{i}^{\ast}(t)=0$ for $t>t_{2}%
$, and the Hamiltonian is discontinuous at $t=t_{2}$ (entry point of the state
constrained arc $v_{i}^{\ast}(t)=v_{max}$); see \citet{bryson1975applied}. It
follows from \eqref{eq:model2} that for $t\geq t_{2}$
\begin{equation}
p_{i}^{\ast}(t)=v_{max}~t+r_{i}, \label{eq:cas3g}%
\end{equation}
where $r_{i}$ is the constant of integration that can be computed from the
position of the vehicle at $t=t_{2}^{-}$.

Given certain terminal constraints, it is possible that the state variable
constraint becomes inactive again; see \citet{bryson1975applied}. If this
happens at time $t=t_{3}>t_{2}$ (exit point of the corner) the state variable
constraint becomes inactive again, i.e., $v_{min}<v_{i}(t)<v_{max}$, then the
Hamiltonian and costates are continuous at $t=t_{3}$, i.e., $H(t_{3}%
^{-})=H(t_{3}^{+})$, $\lambda_{i}^{p}(t_{3}^{-})=\lambda_{i}^{p}(t_{3}%
^{+})=g_{i}$, and $\lambda_{i}^{v}(t_{3}^{-})=\lambda_{i}^{v}(t_{3}%
^{+})=-(g_{i}t+h_{i})$, where $g_{i}$ and $h_{i}$ are constants of
integration. Hence
\begin{equation}
-\frac{1}{2}u_{i}^{\ast}(t)=g_{i}(v_{max}-v_{i}(t)). \label{eq:cas3i2}%
\end{equation}

The optimal control input, speed, and position are
\begin{equation}
\label{eq:cas3j}u^{*}_{i}(t) = g_{i}t +h_{i},
\end{equation}
\begin{equation}
v^{*}_{i}(t) = \frac{1}{2} g_{i}t^{2} + h_{i}t + q_{i} , \label{eq:cas3k}%
\end{equation}
\begin{equation}
p^{*}_{i}(t) = \frac{1}{6} g_{i}t^{3} + \frac{1}{2} h_{i}t^{2} + q_{i}t +
s_{i}, \label{eq:cas3l}%
\end{equation}
where the constants of integration $g_{i}$, $h_{i}$, $q_{i}$, and $s_{i}$ can
be computed from the control, speed, and position of the vehicle at
$t=t_{3}^{-}$ and \eqref{eq:cas3i2} at $t=t_{3}^{+}$.


\textbf{4. Control Constraints Active,} $u_{i}(t)=u_{min}$. Suppose that at
time $t=t_{1}$, (\ref{eq:20}) becomes $u_{i}^{\ast}(t)=u_{min}$ while
$v_{min}<v_{i}(t)<v_{max}$. In this case, the Hamiltonian is continuous at
$t=t_{1}$ (entry point of the control constrained arc). It follows from
\eqref{eq:model2} that for $t\geq t_{1}$
\begin{gather}
v_{i}^{\ast}(t)=u_{min}~t+f_{i},\label{eq:cas4e}\\
p_{i}^{\ast}(t)=\frac{1}{2}u_{min}~t^{2}+f_{i}~t+e_{i}, \label{eq:cas4f}%
\end{gather}
where $f_{i}$ and $e_{i}$ are constants of integration that can be computed
easily since we know the speed and position of the vehicle at time $t=t_{1}$.


\textbf{5. Control and State Constraints Active,} $u_{i}(t)=u_{min}$ and
$v_{i}(t)=v_{min}$. Suppose that at time $t=t_{2}>t_{1}$ (exit point of the
control constrained arc and entry point of the state variable constrained
arc), \eqref{eq:cas4e} becomes equal to $v_{min}$. Then from \eqref{eq:model2}
we have $\dot{v}_{i}^{*}= u_{i}^{*}(t)=0$ for $t>t_{2}$, and the Hamiltonian
is discontinuous at $t=t_{2}$ (corner). Substituting $u_{i}^{*}(t)=0$ for
$t>t_{2}$ into the vehicle dynamics equations \eqref{eq:model2} we can find
the optimal speed and position of each vehicle for $t\ge t_{2}$, namely
\begin{equation}
p_{i}^{\ast}(t)=v_{min}~t+r_{i} \label{eq:cas3g}%
\end{equation}
where $r_{i}$ is the constant of integration that can be computed from the
position of the vehicle at $t=t_{2}^{-}$.

If at time $t=t_{3}>t_{2}$ (exit point of the corner) the state variable
constraint becomes inactive again, i.e., $v_{min} < v_{i}(t) < v_{max}$, then
the Hamiltonian and costates are continuous at $t=t_{3}$. The analysis follows
the discussion at the exit point of the corner in the case where
$u_{i}(t)=u_{max}$ and $v_{i}(t)=v_{max}$, and the optimal control input,
speed, and position are given by \eqref{eq:cas3j}-\eqref{eq:cas3l}.


\textbf{6. State Constraints Active,} $v_{i}(t)=v_{max}$. Suppose that at time
$t=t_{1}$, (\ref{eq:21}) becomes $v_{i}^{\ast}(t)=v_{max}$ while
$u_{min}<u_{i}(t)<u_{max}$. Then from \eqref{eq:model2} we have $\dot{v}%
_{i}^{\ast}=u_{i}^{\ast}(t)=0$ for $t>t_{1}$, and the Hamiltonian is
discontinuous at $t=t_{1}$. Substituting $u_{i}^{\ast}(t)=0$ into the vehicle
dynamics equations \eqref{eq:model2} we can also find the optimal position of
each vehicle for $t\geq t_{1}$, namely
\begin{equation}
p_{i}^{\ast}(t)=v_{max}~t+r_{i} \label{eq:cas3g}%
\end{equation}
where $r_{i}$ is the constant of integration that can be computed from the
position of the vehicle at $t=t_{1}^{-}$.

If at time $t=t_{3}>t_{2}$ (exit point of the state constrained arc) the state
variable constraint becomes inactive again, i.e., $v_{min} < v_{i}(t) <
v_{max}$, then the Hamiltonian and costates are continuous at $t=t_{3}$. The
analysis follows the discussion at the exit point of the state constrained arc
in the case where $u_{i}(t)=u_{max}$ and $v_{i}(t)=v_{max}$, and the optimal
control input, speed, and position are given by \eqref{eq:cas3j}-\eqref{eq:cas3l}.


\textbf{7. State Constraints Active,} $v_{i}(t)=v_{min}$. Suppose that at time
$t=t_{1}$, (\ref{eq:21}) becomes $v_{i}^{\ast}(t)=v_{min}$ (entry point of the
state variable constrained arc) while $u_{min}<u_{i}(t)<u_{max}$. Then from
\eqref{eq:model2} we have $\dot{v}_{i}^{\ast}=u_{i}^{\ast}(t)=0$ for $t>t_{1}%
$, and the Hamiltonian is discontinuous at $t=t_{1}$. It follows from
\eqref{eq:model2} that for $t\geq t_{1}$
\begin{equation}
p_{i}^{\ast}(t)=v_{min}~t+r_{i} \label{eq:cas3g}%
\end{equation}
where $r_{i}$ is the constant of integration that can be computed from the
position of the vehicle at $t=t_{1}^{-}$. The analysis is similar to the case
where $u_{i}(t)=u_{max}$ and $v_{i}(t)=v_{max}$, and the optimal control
input, speed, and position are given by \eqref{eq:cas3j}-\eqref{eq:cas3l}.

To derive the analytical solution of \eqref{eq:decentral}, we follow the standard methodology used in optimal control problems with interior point state and/or control constraints. Namely, we first start with the unconstrained arc and derive the solution using \eqref{eq:23}. If the solution violates any of the state or control constraints, then the unconstrained arc is pieced together with the arc corresponding to the violated constraint, and we re-solve the problem with the two arcs pieced together. The two arcs yield a set of algebraic equations which are solved simultaneously using the boundary conditions of \eqref{eq:decentral} and interior conditions between the arcs. If the resulting solution, which includes the determination of the optimal switching time from one arc to the next one, violates another constraint, then the last two arcs are pieced together with the arc corresponding to the new violated constraint, and we re-solve the problem with the three arcs pieced together. The three arcs will yield a new set of algebraic equations that need to be solved simultaneously using the boundary conditions of \eqref{eq:decentral} and interior conditions between the arcs. The resulting solution includes the optimal switching time from one arc to the next one. The process is repeated until the solution does not violate any other constraints. 


\begin{customrem}{3}
The simple nature of the optimal control and states in (\ref{eq:20}) through
(\ref{eq:22}) makes the online solution of (\ref{eq:decentral})
computationally feasible, even with the additional burden of checking for
active constraints in Cases 2) through 7). However, there is an additional
feature of the solution that we can exploit, i.e., the fact that the control
structure for CAV $i$ remains unchanged until an \textquotedblleft
event\textquotedblright\ e.g., unexpected braking by the
preceding vehicle, rescheduling of the crossing order by the coordinator,
etc.) occurs that affects its behavior. Therefore, there is no need for a
time-driven controller implementation such that $u_{i}^{\ast}(t)$ is
repeatedly re-evaluated. Rather, an event-driven controller may be used
without affecting its optimality properties under conditions such as those
described in \citet{Zhong2010}.
\end{customrem}


\subsection{Feasibility analysis for safety constraints}

\label{sec:4c}

\begin{figure}[ptb]
\centering
\includegraphics[width=3 in]{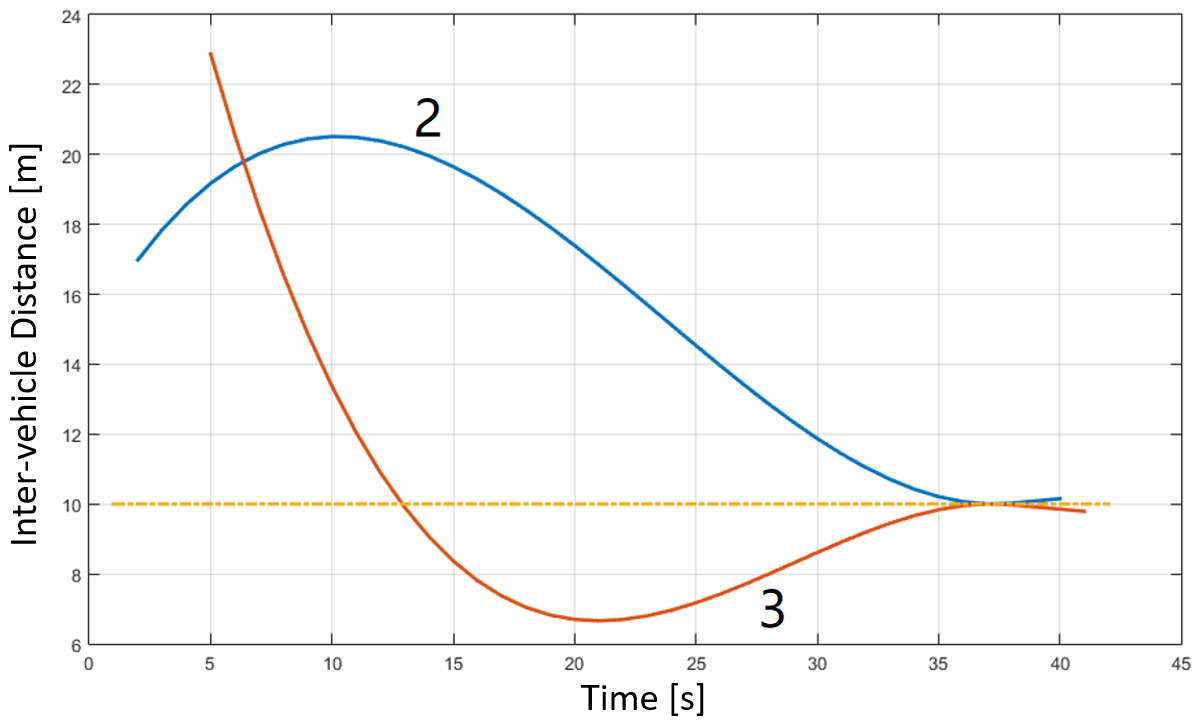}\caption{Example of safety
constraint violation by CAV 3 when $\delta$ = 10.}%
\label{viocase}%
\end{figure}

As already pointed out, the decentralized problem (\ref{eq:decentral}) does
not explicitly include the safety constraints
\eqref{eq:rearend}~and~\eqref{eq:lateral}. While the latter holds by the
construction of $t_{i}^{m}$ in (\ref{def:tf}) and is needed only over
$[t_{i}^{m},t_{i}^{f}]$, the former is not guaranteed to hold for all
$t\in\lbrack t_{i}^{0},t_{i}^{m}]$. We begin with a simple example of how
\eqref{eq:rearend} may be violated under the optimal control (\ref{eq:20}).
This is illustrated in Fig. \ref{viocase} with $\delta=10$ for two CAVs that
follow each other in the same lane within the CZ. We can see that while
\eqref{eq:rearend} is eventually satisfied, due to the constraints imposed on
the solution of \eqref{eq:decentral} through \eqref{eq:lateral}, the
controller (\ref{eq:20}) is unable to maintain \eqref{eq:rearend} throughout
the CZ. What is noteworthy in Fig. \ref{viocase} is that \eqref{eq:rearend} is
violated by CAV 3 at an interval which is \emph{interior} to $[t_{3}^{0}%
,t_{3}^{f}]$, i.e., the form of the optimal control solution (\ref{eq:20})
causes this violation even though the constraint is initially satisfied at
$t_{3}^{0}=5$ in Fig. \ref{viocase}.

Recall that we use $k$ to denote the CAV physically preceding $i$ on the same
lane, and that $i-1$ is the CAV ahead of $i$ in the queue. Clearly,
$k\leqslant i-1$; when $k=i-1$, then $i$ follows $i-1$ in the same lane,
whereas if $i-1$ is on a different lane from $i$, then $k<i-1$. Using this
notation, the following theorem asserts that we can always find initial
conditions $(t_{i}^{0},v_{i}^{0})$ which guarantee the safety constraint
\eqref{eq:rearend} holds throughout the CZ under the decentralized optimal
control, even though \eqref{eq:rearend} is not explicitly included in
(\ref{eq:decentral}).


\begin{customthm}{2}
There exists a nonempty feasible region $\mathcal{F}_{i}\subset R^{2}$ of
initial conditions $(t_{i}^{0},v_{i}^{0})$ for CAV $i$ such that, under the
decentralized optimal control, $s_{i}(t)\geqslant\delta$ holds for all
$t\in\lbrack t_{i}^{0},t_{i}^{m}]$ given initial and final conditions
$t_{k}^{0},v_{k}^{0},t_{k}^{m},v_{k}^{m}$ for CAV $k$. \label{theo:3}
\end{customthm}

\textbf{Proof:} See Appendix.

For any set of initial conditions which are feasible, our analysis gives an
optimal control solution, possibly with a constrained arc. The case which
applies depends on the choice of initial conditions. In other words, our
analysis provides a map from the feasible region to a set of optimal controls
for CAV $i$ which all satisfy the safety inequality. Theorem \ref{theo:3}
asserts that as long as we can drive the CAV to a feasible initial point,
there exists a solution satisfying the safety inequality over the entire CZ
and MZ which may or may not include a constrained arc. There are
two possible ways to deal with the feasibility issue. One approach is to guide
the CAV through an appropriately designed \textquotedblleft Feasibility
Enforcement Zone\textquotedblright\ (FEZ) that precedes the CZ and to make
adjustments so as to attain a feasible initial condition when it reaches the
CZ. The associated feasibility enforcement analysis and the design process of
a FEZ are extensively discussed in \citet{Zhang2016}. Alternatively, if a
FEZ is not realizable and a CAV arrives with $(t_{i}^{0},v_{i}^{0}%
)\notin\mathcal{F}_{i}$, then the decentralized nature of (\ref{eq:decentral})
allows us to forego its optimal control and settle for a non-optimal but safe
control instead with some $t_{i}^{m}$ which is supplied to CAV $i+1$ so as to
continue the use of (\ref{eq:20}) for all subsequent CAVs.


\section{Simulation Results}

To evaluate the effectiveness of the proposed solution, we considered the
following two case studies: (1) coordination of 20 vehicles, (2) coordination
of 448 vehicles. For the first study we used MATLAB and for the second one we
used VISSIM, a microscopic multi-modal traffic flow commercial simulation
software package. The proposed solution was compared to a baseline scenario,
where the intersection has traffic lights with fixed switching times. To
quantify the impact of the vehicle coordination on fuel consumption, we used
the polynomial metamodel proposed in \citet{Kamal2013a} that yields vehicle
fuel consumption as a function of the speed, $v(t)$, and control input, $u(t)$.

\begin{figure}[ptb]
\centering
\includegraphics[width=3 in]{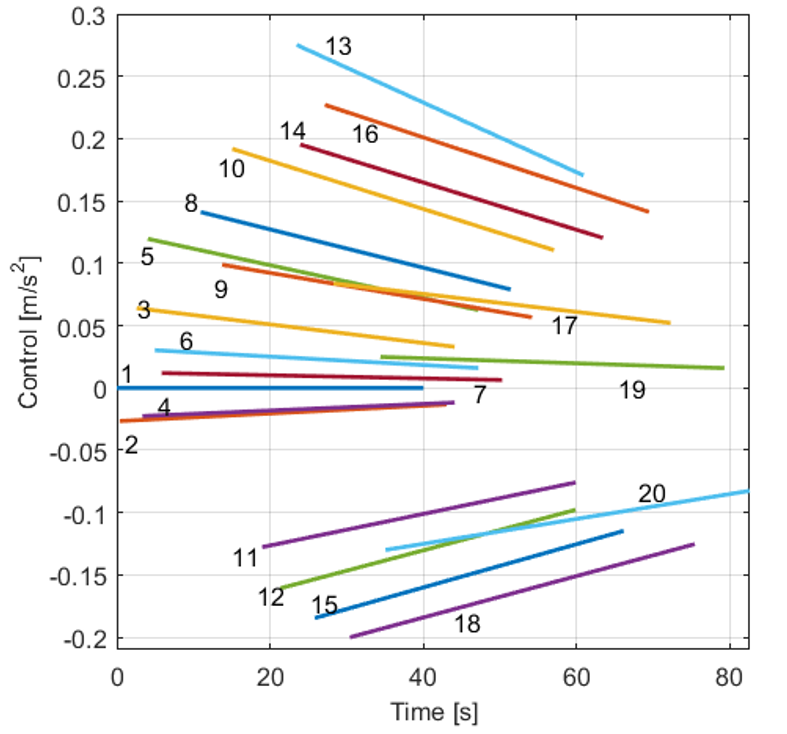} \caption{Optimal control input
signal of the first 20 vehicles in the queue.}%
\label{sim:2}%
\end{figure}

\begin{figure}[ptb]
\centering
\includegraphics[width=3 in]{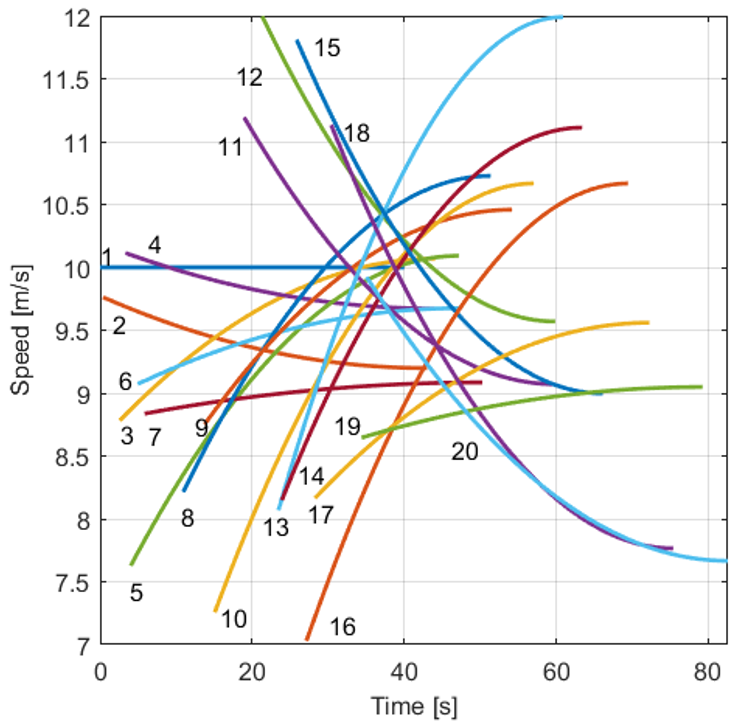} \caption{Speed of the first 20
vehicles in the queue.}%
\label{sim:3}%
\end{figure}

In the first case study, we considered a single intersection, where the length
of the MZ, $S$, is 30m and the length of the CZ, $L$, is 400m. The minimum
safe distance, $\delta,$ between two vehicles was set to be 10m. The maximum
and minimum speed limits are 13 m/s and 0, respectively. The maximum
acceleration limit is 0.2 m/s$^{2}$ and  the maximum deceleration
is set to be arbitrarily large. The control input and the optimal speed for
each vehicle in the queue is shown in Fig. \ref{sim:2} and \ref{sim:3}. Note
that CAV \# 16 violates both the control constraint $u(t)\leq u_{max}$ and the
state constraint $v(t)\leq v_{max}$. 



\begin{figure}[ptb]
\centering
\includegraphics[width=\columnwidth]{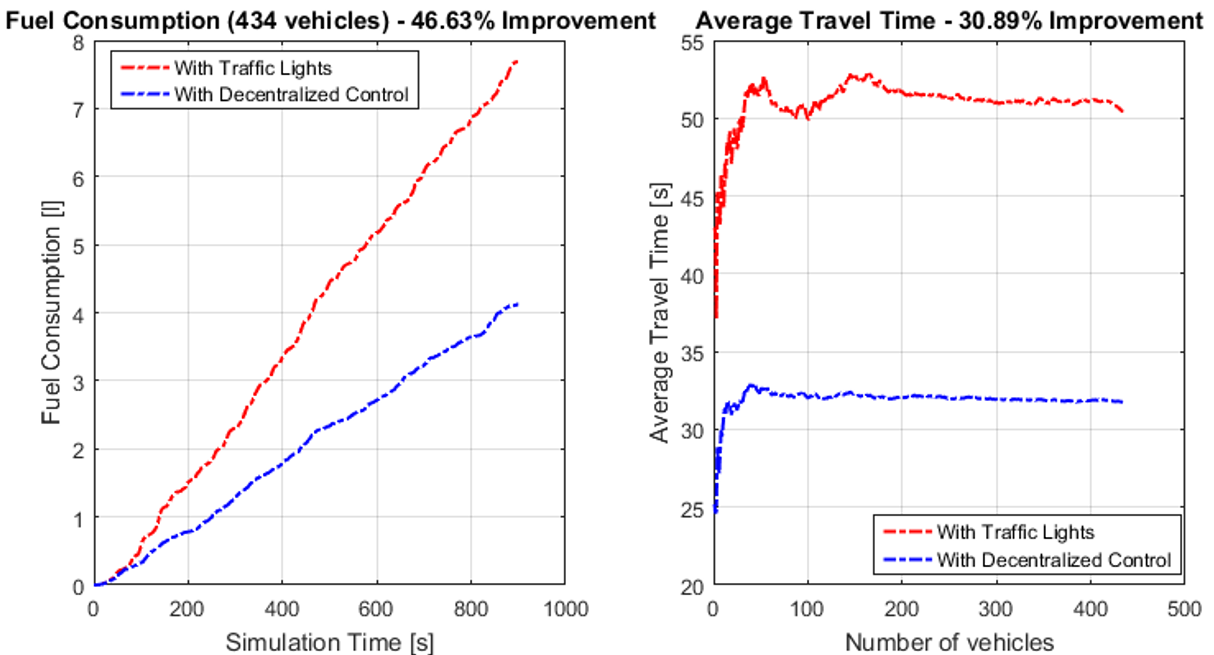} \caption{Fuel
consumption and average travel time improvement.}%
\label{sim:4}%
\end{figure}

In the second case study, we considered two actual intersections in tandem
located in Boston. For each direction, only one lane is considered. We set $L$
= 245 $m$ and $S$ = 35 $m$ for both intersections. As the shapes of the actual
intersections are not regular, the distance between them is not the same for
different directions; in particular, the distance in the lane where the
traffic flow coming from the east is 160 $m$, whereas the distance in the lane
where the traffic flow goes from the west is 145 $m$. In this study, we do not
consider the coupling of the two intersections. The vehicle arrival rate is
assumed to be given by a Poisson process with $\lambda$ = 450 $veh/h$ for each
lane. A comparison to the baseline scenario using traffic lights is shown in
Fig. \ref{sim:4}. The fuel consumption improvement was 46.6\%, while the
travel time was improved by 30.9\%. The fuel consumption improvement is due to
the following reasons: (1) the vehicles do not come to a full stop, thereby
conserving momentum, and (2) each vehicle travels with the minimum
acceleration/deceleration inside the CZ so that transient engine operation is
minimized with direct benefits in fuel consumption.

\section{Concluding Remarks and Future Work}

We have addressed the problem of optimally controlling CAVs crossing an urban
intersection without any explicit traffic signaling. The objective was to
minimize energy consumption subject to a throughput maximization requirement
and hard safety constraints. We have shown that the solution of the latter
depends only on the hard safety constraints imposed on CAVs and that its
structure enables a decentralized energy minimization optimal control problem
formulation whose terminal time depends only on a \textquotedblleft
neighboring\textquotedblright\ set of CAVs. We presented a complete analytical
solution of these decentralized problems and derived conditions under which
feasible solutions satisfying all safety constraints always exist. The
effectiveness of the proposed solution was validated through simulation which
showed that the benefits of the proposed framework are substantial.

In our decentralized framework, we considered full penetration of identical
CAVs having access to perfect information (no errors or delays). We also did
not consider lane changing, turns or pedestrians. Ongoing
research is considering turns [see \citet{Zhang2017}] and lane changing in
the intersection with a diverse set of CAVs and exploring the associated
tradeoffs between the intersection throughput and fuel consumption of each
individual vehicle. Another issue that we are considering is the
potential rear-end collision that would occur inside the MZ if the terminal
speeds of two vehicles  $i$ and $k$ traveling on the same lane are different.
If this case arises, there are two possible approaches to adopt. The first
approach is to set $v_{i}(t_{i}^{m})=v_{k}(t_{k}^{m})$ and solve the optimal
control problem for CAV $i$ with a specified terminal speed. Alternatively, we
can simply forgo the assumption of constant speed in the MZ and ensure the
absence of rear-end collision. We are also investigating the implications of
the proposed approach to adjacent intersections and a feasibility enforcement
to ensure that each CAV starts from a feasible state; see \citet{Zhang2016}.
The fact that the control structure for each CAV remains unchanged until an
\textquotedblleft event" occurs that affects its behavior is an additional
feature of the solution that is being exploited and which will eventually lead
to event-driven controllers. 

The first-in-first-out queue imposes several limitations that can become even more apparent in heavy-volume traffic situations.
An important direction for future research is to relax the first-in-first-out queue and establish a higher-level dynamic optimization problem the solution of which would yield online the optimal ``scheduling" for the CAVs to cross the intersection. Future research should also consider different
penetrations of CAVs, which can alter significantly the efficiency of the
entire system, e.g., what is the critical traffic flow rate beyond which the
benefits of CAVs reach their limit? In fact, as the arrival rates increase,
the optimization process may result in occasional stopping and restarting due
to the implicit state constraint $v_{i}(t) \geq0$. Although it is relatively
straightforward to extend our results to the case where the perfect
information assumption is relaxed, future research needs to be directed at the
implications of errors and/or delays. 

\textbf{Appendix}

\textbf{A. Proof of Lemma 1}

Using the second form of the objective function in (\ref{eq:central1}) and
recalling that $t_{1}^{m^{\ast}}$ is fixed by the unconstrained control
$u_{1}(t)\in\mathcal{U}_{1}$, the solution $\mathbf{t}^{\ast}$ is obtained by
minimizing $t_{i}^{m\ast}$ for every $i=2,\ldots,N$. There are three cases to consider:

\emph{Case 1}: If $i-1\in\mathcal{R}_{i}(t)\cup\mathcal{O}_{i}(t)$, then
\eqref{eq:lateral} cannot become active. On the other hand, \eqref{eq:rearend}
may become active if there exists some $k\in\mathcal{L}_{i}(t)$ such that
$k\in\mathcal{R}_{i-1}(t)\cup\mathcal{O}_{i-1}(t)$, thus creating an
interdependence between $i$ and $i-1$ through $k$. If this happens,
\eqref{eq:rearend} implies that $p_{k}(t_{i}^{m})-p_{i}(t_{i}^{m})=L+v_{k}%
^{m}(t_{i}^{m}-t_{k}^{m}) - L\geqslant\delta$, hence $t_{i}^{m}-t_{k}^{m}%
\geq\frac{\delta}{v_{k}^{m}}$. Therefore, $t_{i}^{m}$ is minimized by setting%
\begin{equation}
t_{i}^{m}=t_{k}^{m}+\frac{\delta}{v_{k}^{m}} \label{case1a}%
\end{equation}
Recalling that $t_{i}^{m}-t_{i-1}^{m}\geq0$ from (\ref{eq:fifo}), it follows
that the optimal value of $t_{i}^{m}$ in this case is given by
\begin{equation}
t_{i}^{m^{\ast}}=\max\{t_{i-1}^{m^{\ast}},t_{k}^{m^{\ast}}+\frac{\delta}%
{v_{k}^{m}}\}. \label{lemmacase1}%
\end{equation}
\emph{Case 2}: If $i-1\in\mathcal{L}_{i}(t)$, then \eqref{eq:lateral} cannot
become active, but constraint \eqref{eq:rearend} can. It follows that
(\ref{lemmacase1}) applies with $k=i-1$, yielding%
\begin{equation}
t_{i}^{m^{\ast}}=t_{i-1}^{m^{\ast}}+\frac{\delta}{v_{i-1}^{m}}.
\label{lemmacase2}%
\end{equation}
\emph{Case 3}: If $i-1\in\mathcal{C}_{i}(t)$, then constraint
\eqref{eq:lateral} can become active. It follows that $t_{i}^{m}\geq
t_{i-1}^{f}=t_{i-1}^{m}+\frac{S}{v_{i-1}^{m}}$. Including the possibility that
\eqref{eq:rearend} becomes active if there exists some $k\in\mathcal{L}%
_{i}(t)$ such that $k\in\mathcal{R}_{i-1}(t)\cup\mathcal{O}_{i-1}(t)$, and
recalling (\ref{eq:fifo}), we have
\begin{gather}
t_{i}^{m^{\ast}} = \max\{t_{i-1}^{m^{\ast}}+\frac{S}{v_{i-1}^{m}}%
,t_{k}^{m^{\ast}}+\frac{\delta}{v_{k}^{m}}\}.
\end{gather}
Observe that if such $k\in\mathcal{L}_{i}(t)$ exists, then $k\in
\mathcal{C}_{i-1}(t)$, since $k$ and $i$ are in the same lane. Therefore,
$i-1$ and $k<i-1$ must also satisfy \eqref{eq:lateral}, i.e., $t_{i-1}%
^{m^{\ast}}\geq t_{k}^{m^{\ast}}+\frac{S}{v_{k}^{m}}$, hence $t_{i-1}%
^{m^{\ast}}+\frac{S}{v_{i-1}^{m}}\geq t_{k}^{m^{\ast}}+\frac{S}{v_{k}^{m}%
}+\frac{S}{v_{i-1}^{m}}>t_{k}^{m^{\ast}}+\frac{\delta}{v_{k}^{m}}$, since
$S>\delta$. It follows that
\begin{equation}
t_{i}^{m^{\ast}}=t_{i-1}^{m^{\ast}}+\frac{S}{v_{i-1}^{m}}. \label{lemmacase3}%
\end{equation}
Combining (\ref{lemmacase1}), (\ref{lemmacase2}) and (\ref{lemmacase3}) we
obtain (\ref{Lemma}). $\square$

\textbf{B. Proof of Theorem 1}

When constraints \eqref{speed_accel constraints} are allowed to be active in
(\ref{eq:central1}), then the values of $t_{i}^{m}$ determined through Lemma
\ref{lemma 2} may not be attainable in minimizing $t_{i}^{m}-t_{1}^{m}$. Thus,
we seek a lower bound to $t_{i}^{m}$, which is independent of these
constraints. There are two cases to consider depending on $t_{i}^{0}$ and on
whether CAV $i$ can reach $v_{max}$ prior to $t_{i-1}^{m}$ or not:

$(i)$ If CAV $i$ enters the CZ at $t_{i}^{0}$, accelerates with $u_{max}$
until it reaches $v_{max}$ and then cruises at this speed until it leaves the
MZ at time $t_{i}^{1f}$, it was shown in
\citet{ZhangMalikopoulosCassandras2016} that $t_{i}^{1f}=t_{i}^{0}+\frac
{L+S}{v_{max}}+\frac{(v_{max}-v_{i}^{0})^{2}}{2u_{max}v_{max}}$. From
Assumption \ref{ass:4}, the time it reaches the MZ is $t_{i}^{1f}-\frac
{S}{v_{max}}=$ $t_{i}^{1}$ in (\ref{ti1def}).

$(ii)$ If CAV $i$ accelerates with $u_{max}$ but reaches the MZ at $t_{i}^{m}$
with speed $v_{i}(t_{i}^{m})<v_{max}$, it was shown in
\citet{ZhangMalikopoulosCassandras2016} that it leaves the MZ at time
$t_{i}^{2f}=t_{i}^{0}+\frac{v_{i}(t_{i}^{m})-v_{i}^{0}}{u_{max}}+\frac
{S}{v_{i}(t_{i}^{m})}$ where $v_{i}(t_{i}^{m})=\sqrt{2Lu_{i,max}+(v_{i}%
^{0})^{2}}$. From Assumption \ref{ass:4}, the time it reaches the MZ is
$t_{i}^{2f}-\frac{S}{v_{i}(t_{i}^{m})}=$ $t_{i}^{2}$ in (\ref{ti2def}).

Thus, $t_{i}^{c}=t_{i}^{1}\mathds{1}_{v_{i}^{m}=v_{max}}+t_{i}^{2}%
(1-\mathds{1}_{v_{i}^{m}=v_{max}})$ ($\mathds{1}$ is the
indicator function) is a lower bound of $t_{i}^{f}$ regardless of the
solution of (\ref{eq:central1}). Combining this lower bound with Lemma
\ref{lemma 2}, we immediately obtain (\ref{def:tf}) including $t_{1}^{m^{\ast
}}$ which is a free variable dependent on $\mathscr{A}_{1}$. $\square$

\textbf{C. Proof of Theorem 2}

To prove the existence of the feasible region, there are two cases to
consider, depending on whether any state or control constraint for either CAV
$i$ or $k$ becomes active in the CZ.

\emph{Case 1:} No state or control constraint is active for either $k$ or $i$
over $[t_{i}^{0},t_{i}^{m}]$. By using (\ref{eq:22}), \eqref{eq:23} and the
definition $s_{i}(t)=p_{k}(t)-p_{i}(t)$, under optimal control we can write%
\begin{gather}
s_{i}(t; t_{i}^{0}, v_{i}^{0}) = s_{i}(t, t_{i}^{m},v_{i}^{m},t_{k}^{0}%
,v_{k}^{0},t_{k}^{m},v_{k}^{m}; t_{i}^{0}, v_{i}^{0})\nonumber\\
= A(t, t_{i}^{m},v_{i}^{m},t_{k}^{0},v_{k}^{0},t_{k}^{m},v_{k}^{m}; t_{i}^{0},
v_{i}^{0})t^{3}\nonumber\\
+ B(t, t_{i}^{m},v_{i}^{m},t_{k}^{0},v_{k}^{0},t_{k}^{m},v_{k}^{m}; t_{i}^{0},
v_{i}^{0})t^{2}\nonumber\\
+ C(t, t_{i}^{m},v_{i}^{m},t_{k}^{0},v_{k}^{0},t_{k}^{m},v_{k}^{m}; t_{i}^{0},
v_{i}^{0})t\nonumber\\
+ D(t, t_{i}^{m},v_{i}^{m},t_{k}^{0},v_{k}^{0},t_{k}^{m},v_{k}^{m}; t_{i}^{0},
v_{i}^{0}), \label{si(t)}%
\end{gather}
where $A$, $B$, $C$ and $D$ are functions defined over $t\in\lbrack t_{i}%
^{0},t_{i}^{m}]$. Recall that CAV $k$ is cruising in the MZ, so that
(\ref{eq:20}) through (\ref{eq:22}) do not apply for $k$ over $[t_{k}%
^{m},t_{i}^{m}]$ leading to different expressions for $A$, $B$, $C$ and $D$.
Therefore, we consider two further subcases, one for $[t_{i}^{0},t_{k}^{m}]$
and the other for $[t_{k}^{m},t_{i}^{m}]$. For ease of notation, in the sequel
we replace $(t_{i}^{0},v_{i}^{0})$ by $(\tau,\upsilon)$.

\emph{Case 1.1:} $t\in\lbrack t_{i}^{0},t_{k}^{m}]$. In this case,
$s_{i}(t;\tau,\upsilon)$ is a cubic polynomial inheriting the cubic structure
of (\ref{eq:22}). We can solve \eqref{eq:23} for the coefficients $a_{k}$,
$b_{k}$, $c_{k}$, $d_{k}$, $a_{i}$, $b_{i}$, $c_{i}$ and $d_{i}$ using the
initial and final conditions of CAVs $k$ and $i$. Then, denoting $A$, $B$, $C$
and $D$ as $A_{1}(\tau,v)$, $B_{1}(\tau,v)$, $C_{1}(\tau,v)$ and $D_{1}%
(\tau,v)$ for $t\in\lbrack t_{i}^{0},t_{k}^{m}]$, these are explicitly given
by%
\begin{gather}
A_{1}(\tau,\upsilon) = \frac{1}{(t_{k}^{0}-t_{k}^{m})^{3}}(2L + (v_{k}%
^{m}+v_{k}^{0})(t_{k}^{0}-t_{k}^{m}))\nonumber\\
- \frac{1}{(\tau-t_{i}^{m})^{3}}(2L + (v_{i}^{m}+\upsilon)(\tau-t_{i}%
^{m})),\nonumber
\end{gather}
\begin{gather}
B_{1}(\tau,\upsilon) = - \frac{1}{(t_{k}^{0}-t_{k}^{m})^{3}}[3L(t_{k}%
^{0}+t_{k}^{m})\nonumber\\
+(v_{k}^{0}(t_{k}^{0}+2t_{k}^{m}) +v_{k}^{m}(2t_{k}^{0}+t_{k}^{m}))(t_{k}%
^{0}-t_{k}^{m})]\nonumber\\
+ \frac{1}{(\tau-t_{i}^{m})^{3}} [3L(\tau+t_{i}^{m})\nonumber\\
+ (\upsilon(\tau+2t_{i}^{m})+v_{i}^{m}(2\tau+t_{i}^{m}))(\tau-t_{i}%
^{m})],\nonumber
\end{gather}
\begin{gather}
C_{1}(\tau,\upsilon) = \frac{1}{(t_{k}^{0}-t_{k}^{m})^{3}}[6t_{k}^{0}
t_{k}^{m} L + [(v_{k}^{0}((t_{k}^{m})^{2}+2t_{k}^{0} t_{k}^{m})\nonumber\\
+ v_{k}^{m}((t_{k}^{0})^{2}+2t_{k}^{m} t_{k}^{0}))](t_{k}^{0}-t_{k}%
^{m})]\nonumber\\
- \frac{1}{(\tau-t_{i}^{m})^{3}}[6\tau t_{i}^{m} L + [(\upsilon((t_{i}%
^{m})^{2}+2\tau t_{i}^{m})\nonumber\\
+ v_{i}^{m}((\tau)^{2}+2t_{i}^{m} \tau))](\tau-t_{i}^{m})],\nonumber
\end{gather}
\begin{gather}
D_{1}(\tau,\upsilon) = \frac{1}{(t_{k}^{0}-t_{k}^{m})^{3}} [L((t_{k}^{0}%
)^{3}-3(t_{k}^{0})^{2} t_{k}^{m})\nonumber\\
- ( v_{k}^{0} t_{k}^{0}(t_{k}^{m})^{2}+ v_{k}^{m} (t_{k}^{0})^{2} t_{k}%
^{m})(t_{k}^{0}-t_{k}^{m})]\nonumber\\
- \frac{1}{(\tau-t_{i}^{m})^{3}} [ L((\tau)^{3}-3(\tau)^{2} t_{i}%
^{m})\nonumber\\
- (\upsilon\tau(t_{i}^{m})^{2}+v_{i}^{m} (\tau)^{2} t_{i}^{m})(\tau-t_{i}%
^{m})]. \label{ABCDcoff}%
\end{gather}

Note that in (\ref{si(t)}) we write $s_{i}(t;\tau,v)$ (recall that $(t_{i}%
^{0},v_{i}^{0})\equiv(\tau,v)$) to emphasize the dependence of $s_{i}(t)$ on
these initial conditions for CAV $i$, i.e., we give a parametric
characterization of $s_{i}(t)$ through $(\tau,v)$. Aside from $(\tau,v)$, the
function $s_{i}(t)$ also depends on two groups of arguments: $(i)$ $t_{k}^{0}%
$, $v_{k}^{0}$, $t_{k}^{m}$ and $v_{k}^{m}\equiv v_{k}(t_{k}^{f})=v_{k}%
(t_{k}^{m})$ are quantities associated with CAV $k$. Since $k<i$, all
information related to this CAV is available and is fixed throughout
$[t_{i}^{0},t_{i}^{m}]$. $(ii)$ $t_{i}^{m}$ and $v_{i}^{m}\equiv v_{i}%
(t_{i}^{f})=v_{i}(t_{i}^{m})$ are quantities which can also be determined
through CAV $k$ or $i-1$.

To summarize, $s_{i}(t;\tau,v)$ varies only with $t$ and $(\tau,v)$ with all
remaining arguments being known to CAV $i$. First, observing that the first
half of each of the coefficient expressions in \eqref{ABCDcoff} (which is
derived by solving \eqref{eq:22} and \eqref{eq:23} for CAV $k$) is a constant
fully determined by information provided by CAV $k$, we can rewrite these as
$K_{A_{1}}$, $K_{B_{1}}$, $K_{C_{1}}$, $K_{D_{1}}$. Therefore, $p_{k}^{\ast
}(t)$ in \eqref{eq:22} can be expressed as
\begin{equation}
p_{k}^{\ast}(t)=K_{A_{1}}t^{3}+K_{B_{1}}t^{2}+K_{C_{1}}t+K_{D_{1}}.
\label{pkt}%
\end{equation}
Next, the second half of the coefficients can be expressed through polynomials
in either $\tau$ or $\upsilon$ explicitly derived by solving \eqref{eq:22} and
\eqref{eq:23} for CAV $i$. We will use the notation $P_{X,n}(\tau)$,
$P_{X,n}(\upsilon)$ to represent polynomials of degree $n=1,2,3$ and
$X\in\{A_{1},B_{1},C_{1},D_{1}\}$. Similarly, we set $Q_{3}(\tau)=(\tau
-t_{i}^{m})^{3}$. Thus, for the coefficients in Eq. \eqref{ABCDcoff}, we get%
\begin{equation}
A_{1}(\tau, \upsilon) = K_{A_{1}} + \frac{ P_{A_{1},1}(\tau) P_{A_{1}%
,1}(\upsilon)}{Q_{3}(\tau)},\nonumber
\end{equation}
\begin{equation}
B_{1}(\tau, \upsilon) = K_{B_{1}} + \frac{ P_{B_{1},2}(\tau)P_{B_{1}%
,1}(\upsilon)}{Q_{3}(\tau)},\nonumber
\end{equation}
\begin{equation}
C_{1}(\tau, \upsilon) = K_{C_{1}} + \frac{ P_{C_{1},3}(\tau) + P_{C_{1}%
,2}(\tau)P_{C_{1},1}(\upsilon)}{Q_{3}(\tau)},\nonumber
\end{equation}
\begin{equation}
D_{1}(\tau, \upsilon) = K_{D_{1}} + \frac{ P_{D_{1},3}(\tau) + P_{D_{1}%
,2}(\tau)P_{D_{1},1}(\upsilon)}{Q_{3}(\tau)}. \label{ABCD}%
\end{equation}
Note that $p_{k}^{\ast}(t)$ in (\ref{pkt}) involves only the $K$ terms, while
the analogous cubic polynomial for $p_{i}^{\ast}(t)$ involves only the $P$ and
$Q$ terms.

Our goal is to ensure that $s_{i}(t;\tau,\upsilon)\geqslant\delta$ for all
$t\in\lbrack\tau,t_{k}^{m}]$ (recall that $t_{i}^{0}\equiv\tau$). We can
guarantee this by ensuring that $s_{i}^{\ast}(\tau,\upsilon)\equiv\min
_{t\in\lbrack\tau,t_{k}^{m}]}\{s_{i}(t;\tau,\upsilon)\}\geqslant\delta$. Thus,
we shift our attention to the determination of $s_{i}^{\ast}(\tau,\upsilon)$.
We can obtain expressions for the first and the second derivative of
$s_{i}(t;\tau,\upsilon)$, $\dot{s}_{i}(t;\tau,\upsilon)$ and $\ddot{s}%
_{i}(t;\tau,\upsilon)$ respectively, from \eqref{si(t)}, as follows:%
\begin{gather}
\dot{s}_{i}(t;\tau,\upsilon)=v_{k}(t)-v_{i}(t)=3A_{1}(\tau,\upsilon
)t^{2}+2B_{1}(\tau,\upsilon)t\nonumber\\
+C_{1}(\tau,\upsilon),\label{si_dot}\\
\ddot{s}_{i}(t;\tau,\upsilon)=u_{k}(t)-u_{i}(t)=6A_{1}(\tau,\upsilon
)t+2B_{1}(\tau,\upsilon). \label{si_dotdot}%
\end{gather}
Clearly, we can determine $t_{i}^{\ast}\equiv\arg\min_{t\in\lbrack\tau
,t_{k}^{m}]}\{s_{i}(t;\tau,\upsilon)\}$ as the solution of $\dot{s}_{i}%
(t;\tau,\upsilon)=0$ with $\ddot{s}_{i}(t;\tau,\upsilon)\geqslant0$, unless
$s_{i}^{\ast}(\tau,\upsilon)$ occurs at the boundaries, i.e., $t_{i}^{\ast
}=\tau$ or $t_{i}^{\ast}=t_{k}^{m}$. Thus, there are three cases to consider:

\emph{Case 1.1.A}: $t_{i}^{\ast}=\tau$. In this case,%
\begin{gather}
s_{i}^{\ast}(\tau,\upsilon)=s_{i}(\tau;\tau,\upsilon)\label{casei}\\
=A_{1}(\tau,\upsilon)\tau^{3}+B_{1}(\tau,\upsilon)\tau^{2}+C_{1}(\tau
,\upsilon)\tau+D_{1}(\tau,\upsilon)\geqslant\delta\nonumber
\end{gather}
and we can satisfy $s_{i}(\tau,\upsilon)\geqslant\delta$ for any $\upsilon$ as
long as a feasible $\tau$ is determined. Since at $t=\tau$, we have
$p_{i}(\tau)=0$ and using the definition of $s_{i}(t)=p_{k}(t)-p_{i}(t)$ and
\eqref{pkt}, we get
\[
s_{i}(\tau)=p_{k}^{\ast}(\tau)=K_{A_{1}}\tau^{3}+K_{B_{1}}\tau^{2}+K_{C_{1}%
}\tau+K_{D_{1}}.
\]
Observe that if
\[
p_{k}(\tau)\geqslant\delta
\]
then CAV $i$ enters the CZ at a safe distance from its preceding CAV $k$ and
since $t_{i}^{\ast}=\tau$, we have $s_{i}(t;\tau,\upsilon)\geqslant\delta$ for
all $t\in\lbrack\tau,t_{k}^{m}]$. Thus, it suffices to select $\tau\geqslant
t_{k}^{\delta}$, where $t_{k}^{\delta}$ is the smallest real root of
$p_{k}(\tau)-\delta=0$.

\emph{Case 1.1.B}: $t_{i}^{\ast}=t_{k}^{m}$. In this case,%
\begin{gather}
s_{i}^{\ast}(\tau,\upsilon)=s_{i}(t_{k}^{m};\tau,\upsilon)\label{Case2_tkm}\\
=A_{1}(\tau,\upsilon)(t_{k}^{m})^{3}+B_{1}(\tau,\upsilon)(t_{k}^{m})^{2}%
+C_{1}(\tau,\upsilon)t_{k}^{m}\nonumber\\
+D_{1}(\tau,\upsilon)\geqslant\delta\nonumber
\end{gather}
Thus, the feasibility region $\mathcal{F}_{i}$ is defined by all
$(\tau,\upsilon)$ such that $s_{i}(t_{k}^{m};\tau,\upsilon)-\delta\geq0$ in
the $(\tau,\upsilon)$ space.

\emph{Case 1.1.C}: $t_{i}^{\ast}=t_{1}\in(\tau,t_{k}^{m})$. This case only
arises if the discriminant $\mathcal{D}_{i}(\tau,\upsilon)$ of (\ref{si_dot})
is positive, i.e.,
\begin{equation}
\mathcal{D}_{i}(\tau,\upsilon)=4B_{1}(\tau,\upsilon)^{2}-12A_{1}(\tau
,\upsilon)C_{1}(\tau,\upsilon)>0 \label{t1}%
\end{equation}
and we get%
\begin{equation}
t_{1}=\frac{-2B_{1}(\tau,\upsilon)\pm\sqrt{\mathcal{D}_{i}(\tau,\upsilon)}%
}{6A_{1}(\tau,\upsilon)} \label{t1root}%
\end{equation}
In addition, we must have
\begin{equation}
\tau<t_{1}<t_{k}^{m},\text{ \ \ }\dot{s}_{i}(t_{1};\tau,\upsilon)=0,\text{
\ \ }\ddot{s}_{i}(t_{1};\tau,\upsilon)\geqslant0 \label{t2}%
\end{equation}
Therefore, the feasibility region $\mathcal{F}_{i}$ is defined by all
$(\tau,\upsilon)$ such that%
\begin{equation}
\begin{aligned} s_{i}^{\ast}(\tau,\upsilon) & =s_{i}(t_{1};\tau,\upsilon)\\ =A_1(\tau,\upsilon)(t_{1})^{3}+ & B_1(\tau,\upsilon)(t_{1})^{2}+C_1(\tau,\upsilon )t_{1}\\ + &D_1(\tau,\upsilon)\geqslant\delta \end{aligned} \label{eq:11c}%
\end{equation}
in conjunction with (\ref{t1root})-(\ref{t2}).

\emph{Case 1.2:} $t\in(t_{k}^{m},t_{i}^{m}]$. Over this interval,
$v_{k}(t)=v_{k}^{m}$ by Assumption \ref{ass:4}. Therefore,
\eqref{eq:20}-\eqref{eq:22} no longer apply: \eqref{eq:20} becomes
$u_{k}^{\ast}(t)=0$, \eqref{eq:21} becomes $v_{k}^{\ast}(t)=v_{k}^{m}$ and
\eqref{eq:22} becomes $p_{k}^{\ast}(t)=L+v_{k}^{m}(t-t_{k}^{m})$. Evaluating
$s_{i}(t)=p_{k}(t)-p_{i}(t)$ in this case yields the following coefficients in
(\ref{ABCDcoff}):%

\begin{equation}
\begin{aligned} A_2(\tau,\upsilon) &= - \frac{1}{(\tau-t_i^m)^3}(2L + (v_i^m+\upsilon)(\tau-t_i^m)),\\ B_2(\tau,\upsilon) &= \frac{1}{(\tau-t_i^m)^3} [3L(\tau+t_i^m) + (\upsilon(\tau+2t_i^m) \\ & +v_i^m(2\tau+t_i^m))(\tau-t_i^m)],\\ C_2(\tau,\upsilon) & = v_k^m - \frac{1}{(\tau-t_i^m)^3}[6\tau t_i^m L + [(\upsilon((t_i^m)^2+2\tau t_i^m) \\ & + v_i^m((\tau)^2+2t_i^m \tau))](\tau-t_i^m)],\\ D_2(\tau,\upsilon) & = L - v_k^m t_k^m - \frac{1}{(\tau-t_i^m)^3} [ L((\tau)^3-3(\tau)^2 t_i^m) \\\ & - (\upsilon \tau (t_i^m)^2+v_i^m (\tau)^2 t_i^m)(\tau-t_i^m)]. \end{aligned} \label{MZphase}%
\end{equation}
It follows that $K_{A_{1}},K_{B_{1}},K_{C_{1}}$ and $K_{D_{1}}$ in
\eqref{ABCD} should be modified accordingly, giving $K_{A_{2}}=K_{B_{2}}=0$,
$K_{C_{2}}=v_{k}^{m}$ and $K_{D_{2}}=L-v_{k}^{m}t_{k}^{m}$. Since we are
assuming that no control or state constraints are active for CAV $i$, the
designated final time $t_{i}^{m}$ under optimal control satisfies
(\ref{Lemma}), i.e., $s_{i}(t_{i}^{m}) \geq\delta$. Thus, we only need to
consider the subcase where $s_{i}^{\ast}(\tau,\upsilon)$ occurs in $(t_{k}%
^{m},t_{i}^{m})$ and we have
\[
t_{i}^{\ast}=t_{2},\text{ \ \ \ \ }t_{2}\in(t_{k}^{m},t_{i}^{m}).
\]
Proceeding as in \emph{Case 1.1.C}, the feasibility region $\mathcal{F}_{i}$
is defined by all $(\tau,\upsilon)$ such that%
\begin{gather}
s_{i}^{\ast}(\tau,\upsilon)=s_{i}(t_{2};\tau,\upsilon)\\
=A_{2}(\tau,\upsilon)(t_{2})^{3}+B_{2}(\tau,\upsilon)(t_{2})^{2}+C_{2}%
(\tau,\upsilon)t_{2}\nonumber\\
+D_{2}(\tau,\upsilon)\geqslant\delta\nonumber
\end{gather}
in conjunction with (\ref{t1root})-(\ref{t2}), with $A_{1}$, $B_{1}$, $C_{1}$
and $D_{1}$ replaced by $A_{2}$, $B_{2}$, $C_{2}$ and $D_{2}$, and with
$\tau<t_{1}<t_{k}^{m}$ replaced by $t_{k}^{m}<t_{2}<t_{i}^{m}$.

\emph{Case 2:} At least one of the state and control constraints is active
over $[\tau,t_{i}^{m}]$. As discussed in Section \ref{sec:4b}, there are
several cases to consider when state and/or control constraints are active.
Since one or both CAVs $k$ and $i$ may experience an active constraint, all
different combinations need to be considered. We analyze a few in what follows
since it is clear that the remaining cases are handled in a similar fashion.

\emph{Case 2.1:} $v_{k}^{\ast}(t)=v_{max}$ over an optimal trajectory arc,
while CAV $i$ is unconstrained. In this case, \eqref{eq:20}-\eqref{eq:22} no
longer apply for CAV $k$ and the coefficients in \eqref{ABCDcoff} are affected
similar to \emph{Case 1.2}, except that the fixed speed $v_{k}^{m}$ is now
$v_{max}$.

\begin{figure}[ptb]
\centering
\includegraphics[width= 3 in]{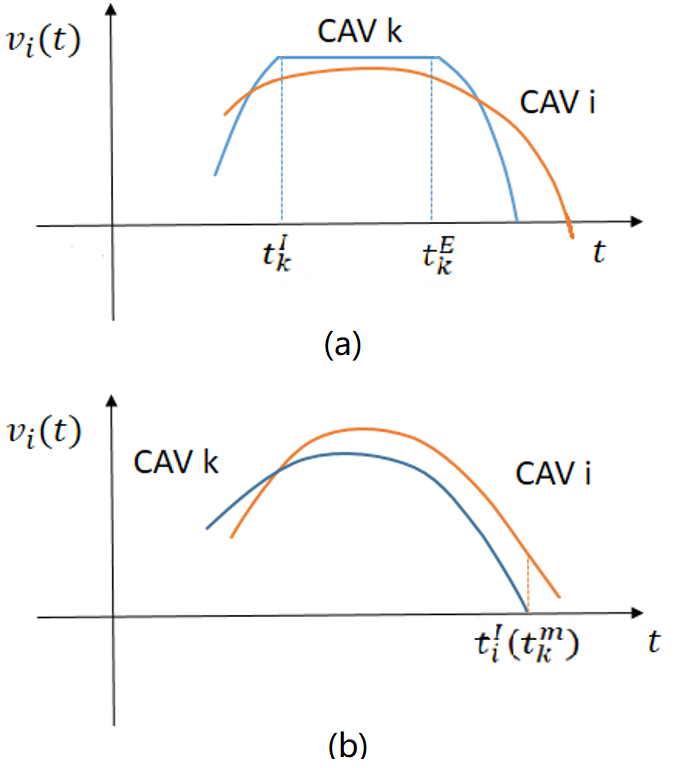}\caption{Cases when the
state and/or control constraints are active.}%
\label{fezcase}%
\end{figure}

First, consider the interval $[\tau,t_{k}^{m}]$. Following the Hamiltonian
analysis in Section \ref{sec:4b}, let $t_{k}^{I}$ be the time CAV $k$ enters
the constrained arc with $v_{k}^{\ast}(t)=v_{max}$ and $t_{k}^{E}$ be the time
it exits this arc (see subfigure (a) in Fig. \ref{fezcase}). The trajectory of
CAV $k$ consists of three arcs as follows. First, for $t\in\lbrack\tau
,t_{k}^{I}]$, $A_{1},B_{1},C_{1}$ and $D_{1}$ are defined exactly as in
\eqref{ABCDcoff}. Second, for $t\in\lbrack t_{k}^{I},t_{k}^{E}]$,
\eqref{eq:20}-\eqref{eq:22} are replaced by $u_{k}^{\ast}(t)=0$, and
$v_{k}^{\ast}(t)=v_{max}$ and \eqref{eq:22} becomes $p_{k}^{\ast}%
(t)=p_{k}(t_{k}^{I})+v_{max}(t-t_{k}^{I})$, where $p_{k}(t_{k}^{I})$ can be
determined before CAV $i$ enters the CZ. The form of the coefficients in
\eqref{ABCDcoff} is modified the same way as in \eqref{MZphase}, with
$v_{k}^{m}$ and $L$ replaced by $v_{max}$ and $p_{k}(t_{k}^{I})$. It follows
that $K_{A_{1}},K_{B_{1}},K_{C_{1}}$ and $K_{D_{1}}$ in \eqref{ABCD} should
also be modified accordingly, with $K_{A_{2}}=K_{B_{2}}=0$, $K_{C_{2}}%
=v_{max}$ and $K_{D_{2}}=p_{k}(t_{k}^{I})-v_{max}t_{k}^{I}$. The final arc is
for $t\in\lbrack t_{k}^{E},t_{k}^{m}]$, when CAV $k$ returns to an
unconstrained arc. The form of the coefficients in \eqref{ABCDcoff} does not
change, except that $A_{1},B_{1},C_{1}$ and $D_{1}$ should be replaced by
$A_{3},B_{3},C_{3}$ and $D_{3}$ since the value of the coefficients may differ
for different unconstrained arcs.

As in \emph{Case 1.1}, we next consider $t_{i}^{\ast}\equiv\arg\min
_{t\in\lbrack\tau,t_{k}^{m}]}$ $\{s_{i}(t;\tau,\upsilon)\}$ and there are
three cases.

\emph{Case 2.1.A:} $t_{i}^{\ast}=\tau$. As in \emph{Case 1.1.A}, it suffices
to select $\tau\geqslant t_{k}^{\delta}$ where $t_{k}^{\delta}$ is the
smallest real root of $p_{k}^{\ast}(\tau)=\delta$.

\emph{Case 2.1.B:} $t_{i}^{\ast}=t_{k}^{m}$. As in \emph{Case 1.1.B}, the
feasibility region $\mathcal{F}_{i}$ is defined by all $(\tau,\upsilon)$ such
that $s_{i}(t_{k}^{m};\tau,\upsilon)-\delta\geq0$ in the $(\tau,\upsilon)$
space, with $A_{1},B_{1},C_{1}$ and $D_{1}$ being replaced by $A_{3}%
,B_{3},C_{3}$ and $D_{3}$.

\emph{Case 2.1.C:} $t_{i}^{\ast}=t_{1}\in(\tau,t_{k}^{m})$. This case may only
arise for $t\in(\tau,t_{k}^{I})$. As in \emph{Case 1.1.C}, the feasibility
region $\mathcal{F}_{i}$ is defined by
\begin{gather}
s_{i}^{\ast}(\tau,\upsilon)=s_{i}(t_{1};\tau,\upsilon)\\
=A_{1}(\tau,\upsilon)(t_{1})^{3}+B_{1}(\tau,\upsilon)(t_{1})^{2}+C_{1}%
(\tau,\upsilon)t_{1}\nonumber\\
+D_{1}(\tau,\upsilon)\geqslant\delta\nonumber
\end{gather}
in conjunction with (\ref{t1root})-(\ref{t2}) with $\tau<t_{1}<t_{k}^{m}$
being replaced by $\tau<t_{1}<t_{k}^{I}$.

For $t\in(t_{k}^{m},t_{i}^{m}]$, the analysis is exactly the same as the way
we handle \emph{Case 1.2}, with $A_{2},B_{2},C_{2}$ and $D_{2}$ being replaced
by $A_{4},B_{4},C_{4}$ and $D_{4}$. The feasibility region $\mathcal{F}_{i}$
is defined by%
\begin{gather}
s_{i}^{\ast}(\tau,\upsilon)=s_{i}(t_{2};\tau,\upsilon)\\
=A_{4}(\tau,\upsilon)(t_{2})^{3}+B_{4}(\tau,\upsilon)(t_{2})^{2}+C_{4}%
(\tau,\upsilon)t_{2}\nonumber\\
+D_{4}(\tau,\upsilon)\geqslant\delta\nonumber
\end{gather}
in conjunction with (\ref{t1root})-(\ref{t2}).

\emph{Case 2.2:} CAV $k$ is unconstrained and $u_{i}^{\ast}(t)=u_{min}$ over
an optimal trajectory arc. Since there are many subcases and each can be
similarly handled, we only consider one subcase where CAV $i$ enters the
constrained arc at $t_{k}^{m}$ (see subfigure (b) in Fig. \ref{fezcase}).

Since both CAV $k$ and $i$ are unconstrained in $[\tau,t_{k}^{m}]$, the form
of the coefficients in \eqref{ABCDcoff} does not change, and the feasibility
region for \emph{Cases 2.2.A, 2.2.B }and\emph{ 2.2.C} can be derived in the
same way as \emph{Case 1.1.A, 1.1.B }and\emph{ 1.1.C}. For $t\in(t_{k}%
^{m},t_{i}^{m}]$, CAV $i$ is deceleraing at a constant value $u_{i}^{\ast
}(t)=u_{min}$. Thus, \eqref{eq:20}-\eqref{eq:22} are replaced by $u_{i}^{\ast
}(t)=u_{min}$, $v_{i}^{\ast}(t)=v_{i}(t_{k}^{m})+u_{min}(t-t_{k}^{m})$ and
$p_{i}^{\ast}(t)=L+v_{i}(t)(t-t_{k}^{m})$, where $v_{i}(t_{k}^{m})$ can be
determined given $(\tau,v)$. The coefficients in \eqref{ABCDcoff} are modifed
as follows:
\begin{equation}
\begin{aligned} A_2(\tau,\upsilon) &= 0, \\ B_2(\tau,\upsilon) &= - \frac{1}{2} u_{min}, \\ C_2(\tau,\upsilon) & = v_k^m - v_i(t_k^m) + u_{min} t_k^m,\\ D_2(\tau,\upsilon) & = - v_k^m t_k^m. \end{aligned}
\end{equation}
For $t\in(t_{k}^{m},t_{i}^{m}]$, CAV $k$ is cruising at the speed $v_{k}^{m}$
and CAV $i$ keeps decelerating until it reaches $v_{i}^{m}$. Therefore,
$s_{i}^{\ast}$ may only occur at $t_{i}^{m}$ and we have
\begin{gather}
s_{i}^{\ast}(\tau,\upsilon)=s_{i}(t_{i}^{m};\tau,\upsilon)\\
=A_{2}(\tau,\upsilon)(t_{i}^{m})^{3}+B_{2}(\tau,\upsilon)(t_{i}^{m})^{2}%
+C_{2}(\tau,\upsilon)t_{i}^{m}\nonumber\\
+D_{2}(\tau,\upsilon)\geqslant\delta\nonumber
\end{gather}
Thus, the feasibility region $\mathcal{F}_{i}$ is defined by all
$(\tau,\upsilon)$ such that \ $s_{i}(t_{i}^{m};\tau,\upsilon)-\delta\geq0$ in
the $(\tau,\upsilon)$ space.

All remaining cases are similarly handled and in each case a feasibility
region $\mathcal{F}_{i}$ is defined by all $(\tau,\upsilon)$ satisfying an
inequality of the form $s_{i}(\sigma;\tau,\upsilon)-\delta\geq0$ for an
appropriate value of $\sigma$ and coefficients in (\ref{si(t)}).

To complete the proof, we show that feasibility region $\mathcal{F}_{i}$ is
always nonempty. This is easily established by considering a point
$(\tau,\upsilon)$ such that $v_{min}<\upsilon<v_{max}$ (which is possible by
Assumption \ref{ass:feas}) and $\tau=t_{k}^{f}$: since $p_{k}^{\ast}(t_{k}%
^{f})=L+S$ and $p_{i}^{\ast}(\tau)=0$, it follows that $s_{i}(\tau)>S>\delta$.
Obviously, any such $(\tau,\upsilon)$ is feasible. $\square$

\begin{figure}[ptb]
\centering
\includegraphics[width=2.8 in]{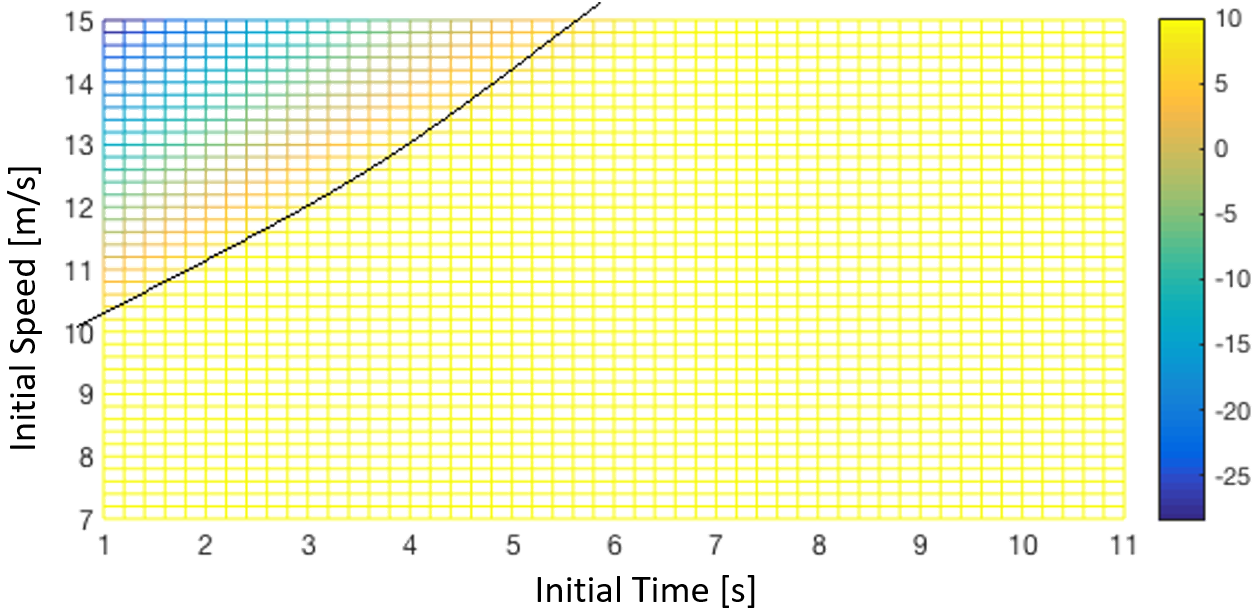}\caption{Feasible and
infeasible region}%
\label{fezregion}%
\end{figure}


\begin{customrem}{4}
To illustrate the feasible region and provide some intuition, we give a
numerical example (see Fig. \ref{fezregion}), with $\delta=10$, $L=400$, and
CAV $k$ is the first vehicle in the CZ and is driving at the constant speed
$v_{k}^{0}=v_{k}^{m} =10$. The color bar in Fig. \ref{fezregion} indicates the
value of $s_{i}^{\ast}(t)$ and the yellow region, determined by
\eqref{eq:11c}, represents the feasible region, while the non-yellow region
represents the infeasible region. The black curve is the boundary between the
two regions and is not linear in general. This boundary curve shifts depending
on the different cases we have considered in the proof of Theorem
\ref{theo:3}. This example also illustrates that we can always find a nonempty
feasible region since we can select points to the right of the curve
corresponding to CAV $i$ entry times in the CZ which can be arbitrarily large.
\end{customrem}


\bibliographystyle{abbrvnat}
\bibliography{TAC_references}

\begin{thebibliography}{34}
\providecommand{\natexlab}[1]{#1}
\providecommand{\url}[1]{\texttt{#1}}
\expandafter\ifx\csname urlstyle\endcsname\relax
  \providecommand{\doi}[1]{doi: #1}\else
  \providecommand{\doi}{doi: \begingroup \urlstyle{rm}\Url}\fi

\bibitem[Alonso et~al.(2011)Alonso, Milan\'{e}s, P\'{e}rez, Onieva,
  Gonz\'{a}lez, and de~Pedro]{Alonso2011}
J.~Alonso, V.~Milan\'{e}s, J.~P\'{e}rez, E.~Onieva, C.~Gonz\'{a}lez, and
  T.~de~Pedro.
\newblock {Autonomous vehicle control systems for safe crossroads}.
\newblock \emph{Transportation Research Part C: Emerging Technologies},
  19\penalty0 (6):\penalty0 1095--1110, Dec. 2011.

\bibitem[Athans(1969)]{Athans1969}
M.~Athans.
\newblock {A unified approach to the vehicle-merging problem}.
\newblock \emph{Transportation Research}, 3\penalty0 (1):\penalty0 123--133,
  1969.
\newblock ISSN 00411647.
\newblock \doi{10.1016/0041-1647(69)90109-9}.

\bibitem[Bryson(1975)]{bryson1975applied}
A.~E. Bryson.
\newblock \emph{Applied optimal control: optimization, estimation and control}.
\newblock CRC Press, 1975.

\bibitem[Colombo and {Del Vecchio}(2014)]{Colombo2014}
A.~Colombo and D.~{Del Vecchio}.
\newblock {Least Restrictive Supervisors for Intersection Collision Avoidance:
  A Scheduling Approach}.
\newblock \emph{IEEE Transactions on Automatic Control}, Provisiona, 2014.

\bibitem[{de La Fortelle}(2010)]{DeLaFortelle2010}
A.~{de La Fortelle}.
\newblock {Analysis of reservation algorithms for cooperative planning at
  intersections}.
\newblock \emph{13th International IEEE Conference on Intelligent
  Transportation Systems}, pages 445--449, Sept. 2010.

\bibitem[Dresner and Stone(2004)]{Dresner2004}
K.~Dresner and P.~Stone.
\newblock {Multiagent traffic management: a reservation-based intersection
  control mechanism}.
\newblock In \emph{Proceedings of the Third International Joint Conference on
  Autonomous Agents and Multiagents Systems}, pages 530--537, 2004.

\bibitem[Dresner and Stone(2008)]{Dresner2008}
K.~Dresner and P.~Stone.
\newblock {A Multiagent Approach to Autonomous Intersection Management}.
\newblock \emph{Journal of Artificial Intelligence Research}, 31:\penalty0
  591--653, 2008.

\bibitem[Huang et~al.(2012)Huang, Sadek, and Zhao]{Huang2012}
S.~Huang, A.~Sadek, and Y.~Zhao.
\newblock {Assessing the Mobility and Environmental Benefits of
  Reservation-Based Intelligent Intersections Using an Integrated Simulator}.
\newblock \emph{IEEE Transactions on Intelligent Transportation Systems},
  13\penalty0 (3):\penalty0 1201,1214, 2012.
\newblock \doi{10.1109/TITS.2012.2186442}.

\bibitem[Kamal et~al.(2013)Kamal, Mukai, Murata, and Kawabe]{Kamal2013a}
M.~Kamal, M.~Mukai, J.~Murata, and T.~Kawabe.
\newblock {Model Predictive Control of Vehicles on Urban Roads for Improved
  Fuel Economy}.
\newblock \emph{IEEE Transactions on Control Systems Technology}, 21\penalty0
  (3):\penalty0 831--841, 2013.
\newblock \doi{10.1109/TCST.2012.2198478}.

\bibitem[Kim and Kumar(2014)]{Kim2014}
K.-D. Kim and P.~Kumar.
\newblock {An MPC-Based Approach to Provable System-Wide Safety and Liveness of
  Autonomous Ground Traffic}.
\newblock \emph{IEEE Transactions on Automatic Control}, 59\penalty0
  (12):\penalty0 3341--3356, 2014.
\newblock \doi{10.1109/TAC.2014.2351911}.

\bibitem[Lee et~al.(2013)Lee, Park, Malakorn, and So]{Lee2013}
J.~Lee, B.~B. Park, K.~Malakorn, and J.~J. So.
\newblock {Sustainability assessments of cooperative vehicle intersection
  control at an urban corridor}.
\newblock \emph{Transportation Research Part C: Emerging Technologies},
  32:\penalty0 193--206, 2013.

\bibitem[Levine and Athans(1966)]{Levine1966}
W.~Levine and M.~Athans.
\newblock {On the optimal error regulation of a string of moving vehicles}.
\newblock \emph{IEEE Transactions on Automatic Control}, 11\penalty0
  (3):\penalty0 355--361, 1966.
\newblock ISSN 0018-9286.
\newblock \doi{10.1109/TAC.1966.1098376}.

\bibitem[Malikopoulos(2011)]{Malikopoulos2011}
A.~A. Malikopoulos.
\newblock \emph{{Real-Time, Self-Learning Identification and Stochastic Optimal
  Control of Advanced Powertrain Systems}}.
\newblock ProQuest, 2011.

\bibitem[Malikopoulos(2013)]{malikopoulos2013stochastic}
A.~A. Malikopoulos.
\newblock Stochastic optimal control for series hybrid electric vehicles.
\newblock In \emph{American Control Conference (ACC), 2013}, pages 1189--1194.
  IEEE, 2013.

\bibitem[Malikopoulos and Aguilar(2013)]{Malikopoulos2013}
A.~A. Malikopoulos and J.~P. Aguilar.
\newblock {An Optimization Framework for Driver Feedback Systems}.
\newblock \emph{IEEE Transactions on Intelligent Transportation Systems},
  14\penalty0 (2):\penalty0 955--964, 2013.
\newblock \doi{10.1109/TITS.2013.2248058}.

\bibitem[Margiotta and Snyder(2011)]{Margiotta2011}
R.~Margiotta and D.~Snyder.
\newblock {An agency guide on how to establish localized congestion mitigation
  programs}.
\newblock Technical report, U.S. Department of Transportation. Federal Highway
  Administration, 2011.

\bibitem[Miculescu and Karaman(2014)]{Miculescu2014}
D.~Miculescu and S.~Karaman.
\newblock {Polling-Systems-Based Control of High-Performance Provably-Safe
  Autonomous Intersections}.
\newblock In \emph{53rd IEEE Conference on Decision and Control}, 2014.

\bibitem[Ntousakis et~al.(2016)Ntousakis, Nikolos, and
  Papageorgiou]{Ntousakis:2016aa}
I.~A. Ntousakis, I.~K. Nikolos, and M.~Papageorgiou.
\newblock Optimal vehicle trajectory planning in the context of cooperative
  merging on highways.
\newblock \emph{Transportation Research Part C: Emerging Technologies},
  71:\penalty0 464--488, 2016.

\bibitem[Rajamani(2012)]{Rajamani2012}
R.~Rajamani.
\newblock \emph{Vehicle Dynamics and Control}.
\newblock Springer, 2012.

\bibitem[Rajamani et~al.(2000)Rajamani, Tan, Law, and Zhang]{Rajamani2000}
R.~Rajamani, H.-S. Tan, B.~K. Law, and W.-B. Zhang.
\newblock {Demonstration of integrated longitudinal and lateral control for the
  operation of automated vehicles in platoons}.
\newblock \emph{IEEE Transactions on Control Systems Technology}, 8\penalty0
  (4):\penalty0 695--708, 2000.

\bibitem[Rios-Torres and Malikopoulos(2017{\natexlab{a}})]{Malikopoulos2016a}
J.~Rios-Torres and A.~A. Malikopoulos.
\newblock A survey on the coordination of connected and automated vehicles at
  intersections and merging at highway on-ramps.
\newblock \emph{IEEE Transactions on Intelligent Transportation Systems},
  18\penalty0 (5):\penalty0 1066--1077, 2017{\natexlab{a}}.

\bibitem[Rios-Torres and Malikopoulos(2017{\natexlab{b}})]{Rios-Torres2}
J.~Rios-Torres and A.~A. Malikopoulos.
\newblock Automated and cooperative vehicle merging at highway on-ramps.
\newblock \emph{IEEE Transactions on Intelligent Transportation Systems},
  18\penalty0 (4):\penalty0 780--789, 2017{\natexlab{b}}.

\bibitem[Rios-Torres et~al.(2015)Rios-Torres, Malikopoulos, and
  Pisu]{Rios-Torres2015}
J.~Rios-Torres, A.~A. Malikopoulos, and P.~Pisu.
\newblock {Online Optimal Control of Connected Vehicles for Efficient Traffic
  Flow at Merging Roads}.
\newblock In \emph{2015 IEEE 18th International Conference on Intelligent
  Transportation Systems}, pages 2432--2437, 2015.

\bibitem[Schrank et~al.(2015)Schrank, Eisele, Lomax, and Bak]{Schrank2015}
B.~Schrank, B.~Eisele, T.~Lomax, and J.~Bak.
\newblock {2015 Urban Mobility Scorecard}.
\newblock Technical report, Texas A\& M Transportation Institute, 2015.

\bibitem[Shladover et~al.(1991)Shladover, Desoer, Hedrick, Tomizuka, Walrand,
  Zhang, McMahon, Peng, Sheikholeslam, and McKeown]{Shladover1991}
S.~E. Shladover, C.~A. Desoer, J.~K. Hedrick, M.~Tomizuka, J.~Walrand, W.-B.
  Zhang, D.~H. McMahon, H.~Peng, S.~Sheikholeslam, and N.~McKeown.
\newblock {Automated vehicle control developments in the PATH program}.
\newblock \emph{IEEE Transactions on Vehicular Technology}, 40\penalty0
  (1):\penalty0 114--130, 1991.

\bibitem[Tachet et~al.(2016)Tachet, Santi, Sobolevsky, Reyes-Castro, Frazzoli,
  Helbing, and Ratti]{Ratti2016}
R.~Tachet, P.~Santi, S.~Sobolevsky, L.~I. Reyes-Castro, E.~Frazzoli,
  D.~Helbing, and C.~Ratti.
\newblock Revisiting street intersections using slot-based systems.
\newblock \emph{PLOS ONE}, 11\penalty0 (3), 2016.

\bibitem[Varaiya(1993)]{Varaiya1993}
P.~Varaiya.
\newblock Smart cars on smart roads: problems of control.
\newblock \emph{IEEE Transactions on Automatic Control}, 38\penalty0
  (2):\penalty0 195--207, 1993.

\bibitem[Yan et~al.(2009)Yan, Dridi, and {El Moudni}]{Yan2009}
F.~Yan, M.~Dridi, and A.~{El Moudni}.
\newblock {Autonomous vehicle sequencing algorithm at isolated intersections}.
\newblock \emph{2009 12th International IEEE Conference on Intelligent
  Transportation Systems}, pages 1--6, 2009.

\bibitem[Zhang et~al.(2016)Zhang, Malikopoulos, and
  Cassandras]{ZhangMalikopoulosCassandras2016}
Y.~Zhang, A.~A. Malikopoulos, and C.~G. Cassandras.
\newblock Optimal control and coordination of connected and automated vehicles
  at urban traffic intersections.
\newblock In \emph{Proceedings of the American Control Conference}, pages
  6227--6232, 2016.

\bibitem[Zhang et~al.(2017{\natexlab{a}})Zhang, Cassandras, and
  Malikopoulos]{Zhang2016}
Y.~Zhang, C.~G. Cassandras, and A.~A. Malikopoulos.
\newblock Optimal control of connected automated vehicles at urban traffic
  intersections: A feasibility enforcement analysis.
\newblock In \emph{Proceedings of the 2017 American Control Conference}, pages
  3548--3553, 2017{\natexlab{a}}.

\bibitem[Zhang et~al.(2017{\natexlab{b}})Zhang, Malikopoulos, and
  Cassandras]{Zhang2017}
Y.~Zhang, A.~A. Malikopoulos, and C.~G. Cassandras.
\newblock Decentralized optimal control for connected automated vehicles at
  intersections including left and right turns.
\newblock In \emph{56th IEEE Conference on Decision and Control}, pages
  4228--4433, 2017{\natexlab{b}}.

\bibitem[Zhong and Cassandras(2010)]{Zhong2010}
M.~Zhong and C.~G. Cassandras.
\newblock Asynchronous distributed optimization with event-driven
  communication.
\newblock \emph{IEEE Transactions on Automatic Control}, 55\penalty0
  (12):\penalty0 2735--2750, 2010.

\bibitem[Zhu and Ukkusuri(2015)]{Zhu2015a}
F.~Zhu and S.~V. Ukkusuri.
\newblock {A linear programming formulation for autonomous intersection control
  within a dynamic traffic assignment and connected vehicle environment}.
\newblock \emph{Transportation Research Part C: Emerging Technologies}, Jan.
  2015.
\newblock ISSN 0968090X.
\newblock \doi{10.1016/j.trc.2015.01.006}.

\bibitem[Zohdy et~al.(2012)Zohdy, Kamalanathsharma, and Rakha]{Zohdy2012}
I.~H. Zohdy, R.~K. Kamalanathsharma, and H.~Rakha.
\newblock {Intersection management for autonomous vehicles using iCACC}.
\newblock \emph{2012 15th International IEEE Conference on Intelligent
  Transportation Systems}, pages 1109--1114, 2012.

\end{thebibliography}


\end{document}